\documentclass[11pt]{article}

\title{A generalized Morse index theorem}
\author{Chaofeng Zhu
$^{1,\;2,\;}$\thanks {Partially supported by FANEDD 200215, 973 Program of MOST, Fok Ying Tung Edu. Funds 91002, LPMC of MOE of
China, and Nankai University.  E-mail: zhucf@nankai.edu.cn }\\ \\
$^{1}$ Nankai Institute of
Mathematics\\
$^{2}$ Key Lab of Pure Mathematics and Combinatorics of Ministry of Education\\  Nankai University, Tianjin 300071\\ The People's Republic of China\\ }

\date{}

\begin{document}

\maketitle

\newtheorem{definition}{Definition}[section]
\newtheorem{theorem}{Theorem}[section]
\newtheorem{proposition}{Proposition}[section]
\newtheorem{notation}{Notation}[section]
\newtheorem{lemma}{Lemma}[section]
\newtheorem{remark}{Remark}[section]
\newtheorem{corollary}{Corollary}[section]
\newtheorem{example}{Example}[section]

\def\s{\section}
\def\ss{\subsection}

\def\d{\begin{definition}}
\def\t{\begin{theorem}}
\def\p{\begin{proposition}}
\def\n{\begin{notation}}
\def\la{\begin{lemma}}
\def\r{\begin{remark}}
\def\c{\begin{corollary}}
\def\ee{\begin{equation}}
\def\aa{\begin{eqnarray}}
\def\y{\begin{eqnarray*}}
\def\bd{\begin{description}}

\newcommand{\proof}[1]{\noindent {\em Proof.}$\quad$ {#1} $\hfill\Box$
                                \vspace{2ex}}
\def\ed{\end{definition}}
\def\et{\end{theorem}}
\def\ep{\end{proposition}}
\def\en{\end{notation}}
\def\el{\end{lemma}}
\def\er{\end{remark}}
\def\ec{\end{corollary}}
\def\eee{\end{equation}}
\def\eaa{\end{eqnarray}}
\def\ey{\end{eqnarray*}}
\def\ebd{\end{description}}

\def\nn{\nonumber}
\def\bp{{\bf Proof.}\hspace{2mm}}
\def\qe{\hfill{\rm Q.E.D.}}
\def\lj{\langle}
\def\rj{\rangle}
\def\dd{\diamond}
\def\ox{\mbox{}}
\def\lb{\label}

\def\K{{\bf K}}
\def\R{{\bf R}}
\def\C{{\bf C}}
\def\Z{{\bf Z}}
\def\N{{\bf N}}
\def\UU{{\bf U}}
\def\DD{{\bf D}}

\def\B{{\cal B}}
\def\Cc{{\cal C}}
\def\CL{{\cal CL}}
\def\Aa{{\cal A}}
\def\Ff{{\cal F}}
\def\Gg{{\cal G}}
\def\Hh{{\cal H}}
\def\Kk{{\cal K}}
\def\Pp{{\cal P}}
\def\Jj{{\cal J}}
\def\Ii{{\cal I}}
\def\Ll{{\cal L}}
\def\Uu{{\cal U}}
\def\g{{\bf g}}

\def\gl{{\rm gl}}
\def\GL{{\rm GL}}
\def\Sp{{\rm Sp}}
\def\sp{{\rm sp}}
\def\U{{\rm U}}
\def\O{{\rm O}}
\def\G{{\rm G}}
\def\H{{\rm H}}
\def\P{{\rm P}}
\def\D{{\rm D}}
\def\T{{\rm T}}
\def\dom{{\rm dom}}
\def\Sa{{\rm S}}
\def\sa{{\rm sa}}
\def\rel{\;{\rm rel.}\;}
\def\vp{\epsilon}
\def\mod{\;{\rm mod}\;}
\def\Lie{{\rm Lie}}
\def\diag{{\rm diag}}
\def\im{{\rm im}\;}
\def\Lag{{\rm Lag}}
\def\Gr{{\rm Gr}}
\def\span{{\rm span}}
\def\Stab{{\rm Stab}\;}
\def\sign{{\rm sign}}
\def\sf{{\rm sf}}
\def\ind{{\rm ind}}
\def\rank{{\rm rank}}
\def\Sg{{\Sp(2n,\C)}}
\def\Na{{\cal N}}
\def\dist{{\rm dist}}
\def\GGG{{\rm Gr}}
\def\Mas{{\rm Mas}}

\begin{abstract}
{\it In this paper, we prove a Morse index theorem for the index
form of even order linear Hamiltonian systems on the closed
interval with reasonable self-adjoint boundary conditions. The
highest order term is assumed to be nondegenerate.}

\end{abstract}

\s{Introduction}\label{s:introduction}

\ss{History}\label{ss:history}

Let $(M,g)$ be an $n$-dimensional Riemannian manifold. The classical Morse Index Theorem states that the number of
conjugate points along a geodesic $\gamma:[a,b]\to M$ counted with multiplicities is equal to the index of the
second variation of the Riemannian action functional $E(c)=\frac{1}{2}\int_a^bg(\dot c,\dot c)dt$ at the critical
point $\gamma$, where $\dot c$ denotes $\frac{d}{dt}c$. Such second variation is called the {\bf index form} for
$E$ at $\gamma$. The theorem has later been extended in several directions (see
\cite{AgSa96,Am61,Du76,PiT00,PiT02,Sm65,Uh73} for versions of this theorem in different contexts). In \cite{Du76}
of 1976, J. J. Duistermaat proved his general Morse index theorem for Lagrangian system with positive definite
second order term and selfadjoint boundary conditions. In \cite{AgSa96} of 1996, A. A. Agrachev and A. V. Sarychev
studied the Morse index and rigidity of the abnormal sub-Riemannian geodesics. In \cite{BeEh79a, BeEh79b} of 1979,
J. K. Beem and P. E. Ehrlich considered the semi-Riemannian case. Later in \cite{He94} of 1994, A. D. Helfer give
a generalization. In \cite{PiT00, PiT02} of 2000, P. Piccione and D. V. Tausk proved a version of the Morse index
theorem for geodesics in semi-Riemannian geodesics with both endpoints varies on two submanifolds of $M$ under
some nondegenerate conditions (cf. \cite[Theorem6.4]{PiT02}). However, such nondegenerate conditions is very
difficult to remove. In \cite{GPiPo03} of 2003, Roberto Giamb{\`o}, Paolo Piccione, Alessandro Portaluri was able
to remove these conditions under the boundary condition of fixed endpoints. Their proof is rather technique and
very difficult to generalize. In \cite{Zh01} of 2001, the author is able to solve these difficulty. However, the
proof is rather technique and hard to follow. It is not clear how the author perturbs a given path of Fredholm
self-adjoint operators to make it with only regular crossings in the degenerate case. In \cite{Ed64} of 1964, the
higher even order case is considered by H. Edwards. He proved a version of Morse index theorem for the even order
Hamiltonian systems on the closed interval with positive definite highest order term and special boundary
condition.

\ss{Set up for regular Lagrangian systems}\label{ss:set up}

Let $M$ be a smooth manifold of dimension $n$, points in its tangent bundle $TM$ will be denoted by $(m,v)$, with
$m\in M$, $v\in T_mM$. Let $f$ be a real-valued $C^3$ function on an open subset $Z$ of $\R\times TM$. Then

\ee \lb{eq05m1.1} E(c)=\int_0^Tf\left(t,c(t),{\dot c}(t)\right)dt
\eee
defines a real-valued $C^2$ function $E$ on
the space of curves

\ee \lb{eq05m1.2} {\cal C}=\left\{c\in C^1([0,T],M);(t,c(t),{\dot c}(t))\in Z \;\mbox{\rm for}\;\mbox{\rm all}\;
t\in[0,T]\right\}. \eee
Equipped with the usual topology of uniform convergence of the curves and their
derivatives, the set ${\cal C}$ has a $C^2$ Banach manifold structure modelled on the Banach space
$C^1([0,T],\R^n)$.

Boundary conditions will be introduced by restriction $E$ to the set of curves
\ee\label{e05m1.3} {\cal
C}_N=\{c\in{\cal C};(c(0),c(T))\in N\},\eee
where $N$ is a given smooth submanifold of $M\times M$. The most
familiar examples are $N=\{m(0),m(T)\}$ and $N=\{(m_1,m_2)\in M\times M;m_1=m_2\}$. In the general case ${\cal
C}_N$ is a smooth submanifold of ${\cal C}$ with its tangential space equal to
\ee\label{e05m1.4} T_c{\cal
C}_N=\left\{\delta c\in C^1([0,T],c^*TM);(\delta c(0),\delta c(T))\in T_{(c(0),c(T))}N\right\}.\eee
$c\in{\cal
C}_N$ is called a {\bf stationary curve} for the boundary condition $N$ if the restriction of $E$ to ${\cal C}_N$
has a stationary point at $c$, i.e., if $DE(c)(\delta c)=0$ for all $\delta c\in T_c{\cal C}_N$. For such a curve
$c$ is of class $C^2$.

Let $c\in{\cal C}_N$ be a stationary curve for the boundary condition $N$. Then the second order differential
$D^2E(c)$ of $E$ at $c$ is symmetric bilinear form on $T_c{\cal C}_N$, which is called the {\bf index form} of $E$
at $c$ with respect to the boundary condition $N$. We want to understand the Morse index of this form, i.e. the
maximal dimension of negative definite subspace of the space $T_c{\cal C}_N$ for the form $D^2E(c)$. In general
the Morse index of the form $D^2E(c)$ on $T_c{\cal C}_N$ will be infinite. In order to get a well-defined integer,
we introduce the following concept.

Assume that $f$ is a {\bf regular Lagrangian}, that is,
\ee\label{e05m1.5} D_v^2f(t,m,v)\; \mbox{is nondegenerate
for all}\; (t,m,v)\in Z.\eee
Here $D_v$ denotes differential of functions on $Z$ with respect to $v\in T_mM$,
keeping $t$ and $m$ fixed. The condition (\ref{e05m1.5}) is called the {\bf Legendre condition}.

Let $H=H^1(T_c{\cal C}_N)$ be the $H^1$ completion of $T_c{\cal C}_N$. By Sobolev embedding theorem, $H\subset
C([0,T],c^*TM)$. Then $D^2E(c)$ is well-defined on $H$. In local coordinates, we have \aa
D^2E(c)(X,Y)&=&\int_0^T\left(D_v^2f(\tilde c(t))({\dot \alpha},{\dot \beta})+D_mD_vf(\tilde
c(t))(\alpha,{\dot \beta})\right.\nn\\
& &+\left.D_vD_mf(\tilde c(t))({\dot \alpha},\beta)+D_m^2f(\tilde c(t))(\alpha,
\beta)\right)dt,\label{e05m1.6}\eaa where $X,Y\in H$, $\alpha$, $\beta$ are the local coordinate expression of
$X$, $Y$ defined by $X=(\alpha,\partial m)$, $Y=(\beta,\partial m)$, $\partial m$ is the natural frame of $T_mM$,
and we use the abbreviation
$$\tilde c(t)=(t,c(t),{\dot c}(t)).$$

In general $\partial m$ and $\alpha$ is not globally well-defined along the curve $c$. Choose a $C^1$ frame $e$ of
$T_c{\cal C}_N$. Such a frame can be obtained by the parallel transformation of the induced connection on $c^*TM$
of a connection on $TM$ (for example, the Levi-Civita connection with respect to the semi-Riemannian metric on
$M$). Then in local coordinates, there is a $C^1$ path $a(t)\in\GL(n,\R)$ such that $\partial m$ at $c(t)$ is the
pairing $(a(t),e(t))$. Note that $a(t)$ is only locally defined in general. Then the vector fields $X,Y\in H$
along $c$ can be written as $X=(x,e)$, $Y=(y,e)$, where $x,y\in H^1([0,T],\R^n)$ and $((x(0),x(T)),(y(0),y(T))\in
R$, $R$ is defined by
$$R=\left\{(x,y)\in\R^{2n};((x,e(0)),(y,e(T)))\in T_{(c(0),c(T))}N\right\}.$$

So we have
\ee \label{e05m1.7}x=a\alpha,\quad {\dot x}=a{\dot\alpha}+{\dot a}\alpha,\quad y=a\beta,\quad {\dot
y}=a{\dot\beta}+{\dot a}\beta.\eee
Substitute (\ref{e05m1.7}) to (\ref{e05m1.6}), we get the following form of the
index form:
\ee \label{e05m1.8} D^2E(c)(X,Y)=\int_0^T\left(\langle p{\dot x}+qx,{\dot y}\rangle+\langle q^*{\dot
x},y\rangle+\langle rx,y\rangle\right)dt,\eee
where $p,q,r\in C([0,T],\gl(n,\R))$, $p$ is of class $C^1$,
$p(t)=p^*(t)$, $r(t)=r^*(t)$, $p(t)$ are invertible for all $t\in[0,T]$, and $*$ denotes the conjugate transpose.

Now define \ee \label{e05m1.9} \Ii_{s,R}(x,y)=\int_0^T\left(\langle p{\dot x}+sqx,{\dot y}\rangle+\langle sq^*{\dot x},y\rangle+\langle
srx,y\rangle\right)dt,\quad s\in[0,1],\eee where $x,y\in H^1([0,T],\R^n)$ and $((x(0),x(T)),((y(0),y(T))\in R$. Since $p$ is of class $C^1$ and $p(t)$ are
nondegenerate, we can associate the path $\Ii_{s,R}$ with a well-defined finite integer, the spectral flow $\sf\{\Ii_{s,R}\}$. Then we can define the relative
Morse index $I(\Ii_{0,R},\Ii_{1,R})$ to be $-\sf\{\Ii_{s,R}\}$. When $p$ is positive definite, $I(\Ii_{0,R},\Ii_{1,R})$ is the Morse index of $D^2E(c)$. Note
that the forms $\Ii_{s,R}$ will depend on the choice of the frame $e$.

\ss{The highlights of the paper}\label{ss:high}

This paper can be viewed as the revised version of \cite{Zh01}. In this paper, we will prove a general version of
Morse index theorem for the index form of even order linear Hamiltonian systems on the closed interval with
reasonable selfadjoint boundary conditions. The highest order term is assumed to be nondegenerate. As a special
case, we prove the Morse index theorem for regular Lagrangian system with selfadjoint boundary conditions. Note
that the index form (see (\ref{e05m1.7}) below) will takes different forms under different choices of the frames
$e$. Then we show how the indices varies under such choices.

Our approach is inspired by the recent papers \cite{BoZh04,BoZh05}
of B. Booss-Bavnbek and the author. We do not use perturbation
method. Our main results can be viewed as pretty much simple
restatement of \cite[Theorem 6.4]{PiT02} and \cite[Theorem
4.9]{GPiPo03} in their cases. Our index theorem does not contain
any assumption on nondegeneracy for the index form. Moreover, we
consider the spectral flow of the paths connected two given index
forms. The index forms in such a path is in general not a compact
perturbation of a given index form. Such phenomenal occurs when we
consider the connected trajectories between two geodesics on the
manifold. These highlights make it easy to apply our index theorem
in the variational problems.

Our paper is arranged as follows. In \S1, we give the background of the problem. In \S2, we state our main results. In \S3, we discuss the properties of the
spectral flow. In \S3, we discuss the properties of the Maslov indices. In \S5, we prove our main results. In this paper, $\dim$ denotes the complex dimension
if no special description.

\s{Main results}\label{s:main}

We shall consider the general case of even order linear Hamiltonian systems. We will consider the complex case.
The real case is a obvious consequence of the complex case.

Let $m,n\in\Z^+$ be positive integers, and $T\in\R^+$ be a
positive real number. Let $p_{k,l}(s,t)\in\gl(n,\C)$,
$(s,t)\in[0,1]\times[0,T]$ be $(m+1)^2$ continuous families of
matrices, where $k,l=0,\ldots,m$. Assume that for all
$(s,t)\in[0,1]\times[0,T]$,
$p_s(t)=(p_{m-k,m-l}(s,t))_{k,l=0,\ldots,m}\in\gl((m+1)n,\C)$ are
selfadjoint, and $p_{m,m}(s,t)$ are nondegenerate. Assume further
that for all $s\in[0,1]$ and $k,l=0\ldots,m$, $p_{k,l}(s,\cdot)\in
C^{\max\{k,l\}}([0,T],\gl(n,\C))$. Then we have a continuous
family of quadratic forms \ee \label{eq05m2.1}
\Ii_s(x,y)=\int_0^T\left(\sum_{k,l=0}^m \langle
p_{k,l}(s,t)\frac{d^l}{dt^l}x,\frac{d^k}{dt^k}y\rangle\right)dt,
\quad\forall x,y\in H^m([0,T];\C^n). \eee Here $\langle
\cdot,\cdot\rangle$ denotes the standard Hermitian inner product
in $\C^n$, and the norm of the Sobolev space $H^m([0,T];\C^n)$ is
defined by
$$\langle x,y\rangle_m=\int_0^T\left(\sum_{k=0}^m \langle
\frac{d^k}{dt^k}x,\frac{d^k}{dt^k}y\rangle\right)dt, \quad\forall x,y\in H^m([0,T];\C^n).$$

Then we define the boundary condition. Let $R\subset\C^{2mn}$ be a given linear subspace. Define
\ee
\label{e05m2.2}H_R=\left\{x\in H^m([0,T];\C^n);(\frac{d^{m-1}}{dt^{m-1}}x(0),\ldots,x(0),
\frac{d^{m-1}}{dt^{m-1}}x(T),\ldots,x(T))\in R\right\}.\eee
Let $\Ii_{s,R}$ be the restriction of $\Ii_s$ to
$H_R$. The central problem in this paper is to understand the Morse index of the form $\Ii_{1,R}$, i.e. the
maximal dimension of negative definite subspace of the form $\Ii_{1,R}$. As in \S\ref{ss:set up}, we shall use the
minus spectral flow $-\sf\{\Ii_{s,R}\}$ as the "difference" between the "Morse indices" of the forms $\Ii_{1,R}$
and $\Ii_{0,R}$.

Let $L_s$ be the unbounded operator on $L^2([0,T];\C^n)$ with
domain $H^{2m}([0,T];\C^n)$ defined by \ee\label{e05m2.3}
(L_sx)(t)=\sum_{k,l=0}^m(-1)^k\frac{d^k}{dt^k}\left(p_{k,l}(s,t)\frac{d^l}{dt^l}x(t)\right),
\quad \forall x\in H^{2m}([0,T];\C^n).\eee Define $R^{2m,b}$ and
$W_{2m}(R)$ by \aa
R^{2m,b}&=&\left\{(x_1,\ldots,x_{2m})\in\C^{2mn};
\sum_{k=1}^m(-1)^{k-1}\langle x_k,y_{m-k+1}\rangle\right.\nn\\
& &\left.+\sum_{k=m+1}^{2m}(-1)^{k-m}\langle x_k,y_{3m-k+1}\rangle=0\;\mbox{for all}\;(y_1,\ldots,y_{2m})\in R
\right\},\label{e05m2.4}\\
W_{2m}(R)&=&\left\{(x_1,x_2,x_3,x_4)\in\C^{4mn};x_1,x_2,x_3,x_4\in\C^{mn},(x_1,x_3)\in R^{2m,b},(x_2,x_4)\in
R\right\}.\label{e05m2.5}\eaa
For each $x\in H^{2m}([0,T];\C^n)$, let $u_{p_s,x}\in H^1([0,T];\C^{2mn})$, $\tilde
u_{p_s,x}$ and $u_{p_s,x}^k$, $k=0,\ldots,2m$ be defined by
\aa
u_{p_s,x}(t)&=&(u_{p_s,x}^{2m-1}(t),\ldots,u_{p_s,x}^0(t)),\nn\\
u_{p_s,x}^k(t)&=&\frac{d^k}{dt^k}x(t),\qquad \qquad \qquad \qquad
\qquad
\qquad \qquad \qquad \qquad \qquad \qquad \;k=0,\ldots,m-1,\nn\\
u_{p_s,x}^k(t)&=&\sum_{2m-k\le \alpha\le m,0\le\beta\le
m}(-1)^{\alpha-m}\frac{d^{\alpha+k-2m}}{dt^{\alpha+k-2m}}\left(p_{\alpha,\beta}(s,t)
\frac{d^{\beta}}{dt^{\beta}}x(t)\right), \quad
k=m,\ldots,2m.\label{e05m2.6}\eaa Let $L_{s,W_{2m}(R)}$ be the
restriction of $L_s$ on the domain
\[\{x\in H^{2m}([0,T];\C^n);(u_{p_s,x}(0),u_{p_s,x}(T))\in W_{2m}(R)\}.\]
By Lemma 3.5 of \cite{BoZh05}, $L_{s,W_{2m}(R)}$, $0\le s\le 1$ is
a continuous family (in the gap norm sense) of unbounded
selfadjoint Fredholm operators. Again we associate the path with
the minus spectral flow $-\sf\{L_{s,W_{2m}(R)}\}$.

Let $J_{2m,n}\in\GL(2mn,\C)$ be the matrix $(j_{k,l})_{k,l=0,\ldots,2m-1}$, where $j_{k,l}=0_n$ for $k+l\ne 2m-1$, $j_{k,l}=(-1)^{k+m}I_n$ for $k+l=2m-1$, and
we denote by $I_n$ and $0_n$ the identity matrix and the zero matrix on $\C^n$ respectively. When there is no confusion, we will omit the subindex $n$ of $I_n$
and $0_n$. Set
$${\bar u}_{p_s,x}=(u^m_{p_s,x},\ldots,u^0_{p_s,x}),\quad {\bar u}_{0,x}=(\frac{d^m}{dt^m}x,\ldots,x).$$
From (\ref{e05m2.6}), we can define the matrices $U(p_s(t))$ and $V(p_s(t))$ for each $(s,t)\in [0,1]\times[0,T]$ by \ee\label{e05m2.7}{\bar
u}_{p_s,x}(t)=U(p_s(t)){\bar u}_{0,x}(t),\quad {\bar u}_{0,x}(t)=V(p_s(t)){\bar u}_{p_s,x}(t).\eee Let $\Theta_{2m,n}\in\gl(2mn,\C)$ be the matrix
$(\theta_{k,l})_{k,l=0,\ldots,2m-1}$, where $\theta_{k,l}=0_n$ for $k+l\ne 2m-2$ or one of $k=l=m-1$, $\theta_{k,l}=(-1)^{k+m+1}I_n$ for $k+l=2m-2$ and $k,l\ne
m-1$. For each $(s,t)\in [0,1]\times[0,T]$, define the matrices $P(p_s(t))$ and $b(p_s(t))$ in $\gl((m+1)n,\C)$ by \aa \label{e05m2.8}
P(p_s(t))&=&(P_{k,l}(s,t))_{k,l=0,\ldots,m}\\
\label{e05m2.9}b(p_s(t))&=&\Theta_{2m,n}+\diag(0_{(m-1)n},P(p_s(t))),\eaa
where \begin{eqnarray*}P_{0,0}(s,t)&=&p_{m,m}(s,t)^{-1},\\
P_{0,l}(s,t)&=&-p_{m,m}(s,t)^{-1}p_{m,m-l}(s,t),\\
P_{k,0}(s,t)&=&-p_{m-k,m}(s,t)p_{m,m}(s,t)^{-1},\\
P_{k,l}(s,t)&=&p_{m-k,m-l}(s,t)-p_{m-k,m}(s,t)p_{m,m}(s,t)^{-1}p_{m,m-l}(s,t)\end{eqnarray*}
for $k,l=1,\ldots,m$, and we denote by $A^*$ the conjugate
transpose of $A$. For each $s\in[0,1]$, let $\gamma_{p_s}(t)$ be
the fundamental solution of the linear Hamiltonian system
\ee\label{e05m2.10} \dot u=J_{2m,n}b(p_s)u.\eee Then
$\gamma_{p_s}(t)$ are symplectic matrices. Then we can associate
the symplectic path $\gamma_{p_s}(t)$, $0\le t\le T$ with the
Maslov-type index $i_{W_{2m}(R)}(\gamma_{p_s})$ for each
$s\in[0,1]$.

We want to address the following problems for even order case in this paper:
\begin{itemize}
\item give the relationship between the integers $-\sf\{\Ii_{s,R}\}$, $-\sf\{L_{s,W_{2m}(R)}\}$ and $i_{W_{2m}(R)}(\gamma_{p_s})$ for $0\le s\le 1$;

\item calculate $i_{W_{2m}(R)}(\gamma_{p_0})$ for $p_0(t)=\diag(p_{0,0}(0,t),0_{mn})$;

\item for two different choices of the frame $e$, the resulted index form $\Ii_{s,R}$ defined by (\ref{e05m1.9}) will have different forms. In this case,
calculate the difference between the resulted integers $i_{W_{2}(R)}(\gamma_{p_1})$.
\end{itemize}

The following three theorems solve the above problems.

\begin{theorem}\label{t05m2.1} Let $\sf\{\Ii_{s,R},0\le s\le 1\}$ be the spectral flow of $\Ii_{s,R}$,
$\sf\{L_{s,W_{2m}(R)},0\le s\le 1\}$ be the spectral flow of $L_{s,W_{2m}(R)}$, and $i_{W_{2m}(R)}(\gamma_{p_s})$ be the Maslov-type index of $\gamma_{p_s}$
defined below. Then we have
\begin{equation}\label{e05m2.11}
-\sf\{\Ii_{s,R},0\le s\le 1\}=-\sf\{L_{s,W_{2m}(R)},0\le s\le 1\}=i_{W_{2m}(R)}(\gamma_{p_1})-i_{W_{2m}(R)}(\gamma_{p_0}).
\end{equation}
\end{theorem}

Assume that $p_0(t)=\diag(p_{m,m}(0,t),0_{mn})$ for all $t\in[0,T]$. Then we have $(P({p_0}))(t)=(p_0(t))^{-1}$, $b(p_0)(t)=(b_{k,l}(t))_{k,l=0,\ldots,2m-1}$,
and $\gamma_{p_0}(t)=(\gamma_{k,l}(t))_{k,l=0,\ldots,2m-1}$, where $b_{k,l}(t)=0_n$ for $k-l\ne 1$, $b_{k,l}(t)=I_n$ for $k-l=1$ and $k\ne m$, $b_{m,m-1}(t)
=(p_{m,m}(0,t))^{-1}$, $\gamma_{k,l}(t)=0$ for $k<l$, $\gamma_{k,l}(t)=\frac{t^{k-l}}{(k-l)!}I_n$ for $k\ge l$ and $k\le m-1$, or $k\ge l$ and $l\ge m$, and
\begin{eqnarray*}
\gamma_{k,l}(t)&=&\frac{1}{(m-l-1)!}\int_0^t dt_{k-m}\int_0^{t_{k-m}}dt_{k-m-1}\ldots\int_0^{t_1}t_0^{m-l-1}(p_{m,m}(0,t_0))^{-1}dt_0\\
&=&\frac{1}{(k-m)!(m-l-1)!}\int_0^t s^{m-l-1}(t-s)^{k-m}(p_{m,m}(0,s))^{-1}ds
\end{eqnarray*}
for $k\ge m$ and $l\le m-1$.

The form of our symplectic path $\gamma_{p_0}(t)$ looks rather complicated. We will consider the following more general situation to simplify our problem.

Let $K\in\GL(n,\C)$. Set $J_K=\pmatrix{0&-K^*\cr\\ K&0\cr}$. Then
$(\C^{2n},\langle J_K\cdot,\cdot\rangle)$ is a symplectic space.
Let
$\gamma(t)=\pmatrix{M_{1,1}(t)&0\cr\\M_{2,1}(t)&M_{2,2}(t)\cr}$,
$0\le t\le T$ be a path in $\GL(2n,\C)$ with
$M_{2,2}(t)^*KM_{1,1}(t)=K$ and $M_{1,1}(t)^*K^*M_{2,1}(t)$
self-adjoint for each $t\in[0,T]$. Then $\gamma(t)$ is a
symplectic path, i.e., $\gamma(t)^*J_K\gamma(t)=J_K$. Let
$R\subset\C^{2mn}$ be a given linear subspace. Define $R^K$ and
$W_{K}(R)$ by
\begin{eqnarray} R^K&=&\left\{(x_1,x_2)\in\C^{2n};
\langle Kx_1,y_1\rangle-\langle Kx_2,y_2\rangle=0\;\mbox{for all}\;(y_1,y_2)\in R
\right\},\label{e05m2.12}\\
W_{K}(R)&=&\left\{(x_1,x_2,x_3,x_4)\in\C^{4n};x_1,x_2,x_3,x_4\in\C^n,(x_1,x_3)\in R^K,(x_2,x_4)\in R\right\}.\label{e05m2.13}
\end{eqnarray}

\begin{theorem}\label{t05m2.2} For the syplectic path $\gamma$ and the Lagrangian space $W_K(R)$ defined above,
we have
\begin{eqnarray}
\dim(\Gr(\gamma(t))\cap
W_K(R))&=&\dim\ker\left((M_{1,1}(T)^*K^*M_{2,1}(t))|_{S(t)}\right)+\dim
S(t)\nn\\
& &+\dim(\Gr(I_{mn})\cap R)-\dim(\Gr(I_{mn})\cap R^K),\label{e05m2.14}\\
i_{W_K(R)}(\gamma)&=&m^+\left((M_{1,1}(T)^*K^*M_{2,1}(T))|_{S(T)}\right)\nn\\
&
&-m^+\left((M_{1,1}(0)^*K^*M_{2,1}(0))|_{S(0)}\right)\nn\\
& &+\dim S(0)-\dim S(T).\label{e05m2.15}
\end{eqnarray}
where $m^+$ denotes the Morse positive index, and
\[S(t)=\{x\in\C^n;(x,M_{1,1}(t)x)\in R^K\}.\]
\end{theorem}

In our case, set $K_{m,n}=(k_{k,l})_{k,l=0,\ldots,m-1}$, where
$k_{k,l}=0_n$ for $k+l\ne m-1$, $k_{k,l}=(-1)^lI_n$ for $k+l=m-1$.
Then we have $R^{K_{m,n}}=R^{2m,b}$ and
$W_{2m}(R)=W_{K_{m,n}}(R)$. Moreover for the symplectic path
$\gamma=\gamma_{p(0)}$, we have
\begin{equation}\label{e05m2.16}
M_{1,1}(T)^*K_{m,n}^*M_{2,1}(T)=\left(\frac{1}{(m-k-1)!(m-l-1)!}\int_0^Tt^{2m-k-l-2}(p_{m,m}(0,t))^{-1}dt\right)_{k,l=0,\ldots,m-1}.
\end{equation}

As a special case, we get the following higher order generalization of theorem of J. J. Duistermaat \cite{Du76}.

\begin{corollary}\label{c05m2.1} Assume that $p_{m,m}(1,t)$ is positive definite for each $t\in[0,T]$.
Set $p_0(t)=\diag(p_{m,m}(1,t),0_{mn})$. Then we have
\begin{equation}\label{e05m2.17}
m^-(\Ii_{1,R})=m^-(L_{1,W_{2m}(R)})=i_{W_{2m}(R)}(\gamma_{p_1})-\dim
S,
\end{equation}
where $m^-$ denotes the Morse (negative) index, and
\begin{eqnarray*}
\gamma_{p_0}(t)&=&{\left(\begin{array}{cc}
  M_{1,1}(t) & 0 \\
  M_{2,1}(t) & M_{2,2}(t)
  \end{array}\right)},\\
S&=&\{x\in\C^{mn};(x,x)\in R^{2m,b}\}.
\end{eqnarray*}
\end{corollary}

Now we consider the third problem. Then $m=1$ and everything is real. Let $a(t)$ be a $C^1$ path in $\GL(n,\R)$, and
\[R^{'}=\{(x,y)\in\R^{2n};(a(0)x,a(T)y)\in R\}.\]
After the change of the frame $e\mapsto a^{*-1}e$, we have $x\mapsto ax$ and the quadratic form $\Ii_{1,R}$ is changed to the restriction of the form
$\Ii_1(ax,ay)$ on $H_{R^{'}}$. Then we get the corresponding $p^{'}$, $q^{'}$ and $r^{'}$. Set $p_1=\pmatrix{p&q\cr\\ q^*&r\cr}$ and
$p_1^{'}=\pmatrix{p^{'}&q^{'}\cr\\(q^{'})^*&r^{'}\cr}$. Let $\gamma_{p_1}$ and $\gamma_{p_1^{'}}$ be defined by (\ref{e05m2.10}). Then we can prove
\begin{equation}\label{e05m2.18}
\gamma_{p_1^{'}}=\diag(a^*,a^{-1})\gamma_{p_1}\diag(a(0)^{*-1},a(0)).
\end{equation}

\begin{theorem}\label{t05m2.3}
Let $a(t)$, $0\le t\le T$ be a path in $\GL(n,\C)$, and
\[R^{'}=\{(x,y)\in\C^{2n};(a(0)x,a(T)y)\in R\}.\]
Let $\gamma$ be a symplectic path, i.e., $\gamma(t)^*J_{2,n}\gamma(t)=J_{2,n}$ for all $0\le t\le T$. Define the symplectic path $\gamma^{'}$ by
\begin{equation}\label{e05m2.19}
\gamma^{'}=\diag(a^*,a^{-1})\gamma_{p_1}\diag(a(0)^{*-1},a(0)).
\end{equation}
Then we have
\begin{equation}\label{e02m2.20}
i_{W_2(R^{'})}(\gamma^{'})-i_{W_2(R)}(\gamma)=\dim (\Gr(I_n)\cap
(R^{'})^{2,b})) -\dim(\Gr(I_n)\cap R^{2,b}).
\end{equation}
\end{theorem}

\s{Spectral flow}\label{s:spectral flow}

\ss{Definition of the spectral flow}\label{ss:def spectral flow}

Roughly speaking, the spectral flow counts the net number of eigenvalues changing from the negative real half axis
to the non-negative one. The definition goes back to a famous paper by M.~Atiyah, V.~Patodi, and I.~Singer
\cite{AtPaSi76}, and was made rigorous by J.~Phillips \cite{Ph96} for continuous paths of bounded self-adjoint
Fredholm operators, by C. Zhu and Y. Long \cite{ZhLo99} in various non-self-adjoint cases, and by B.
Booss-Bavnbek, M. Lesch, and J. Phillips \cite{BoLePh01} in the unbounded self-adjoint case.

Let $X$ be a complex Hilbert space. For a self-adjoint Fredholm operator $A$ on $X$, there exists a unique
orthogonal decomposition
\begin{equation}\label{e05m3.1}
X=X^+(A)\oplus X^0(A)\oplus X^-(A)
\end{equation}
such that $X^+(A)$, $X^0(A)$ and $X^-(A)$ are invariant subspaces associated to $A$, and $A|_{X^+(A)}$,
$A|_{X^0(A)}$ and $A|_{X^-(A)}$ are positive definite, zero and negative definite respectively. We introduce
vanishing, natural, or infinite numbers
$$
m^+(A):=\dim X^+(A),\ m^0(A):=\dim X^0(A),\ m^-(A):=\dim X^-(A),
$$
and call them {\bf Morse positive index}, {\bf nullity} and {\bf Morse index} of $A$ respectively. For finite-dimensional $X$, the {\bf signature} of $A$ is
defined by $\sign(A)=m^+(A)-m^-(A)$ which yields an integer. The {\bf APS projection} $Q_A$ (where APS stands for Atiyah-Patodi-Singer) is defined by
$$Q_A(x^++x^0+x^-):=x^++x^0,$$
for all $x^+\in X^+(A), x_0\in X^0(A), x^-\in X^-(A)$.

Let $\{A_s\}$, $0\le s\le 1$ be a continuous family of self-adjoint Fredholm operators. The spectral flow
$\sf\{A_s\}$ of the family should be equal to $m^-(A_0)-m^-(A_1)$ if $\dim X<+\infty$. We will generalize this
definition to general Banach space $X$ and general continuous family of admissible operators defined below.

Let $X$ be a complex Banach space. We denote the set of closed operators, bounded linear operators and compact
linear operators on $X$ by $\Cc(X)$, $\B(X)$ and $\CL(X)$ respectively. We will denote the spectrum, the regular
set and the domain of an operator $A\in\Cc(X)$ by $\sigma(A)$, $\rho(A)$ and $\dom(A)$ respectively. Let $N$ be an
bounded open subset of $\C$ and $A\in\Cc(X)$. If there exists an bounded open subset $\tilde N\subset N$ with
$C^1$ boundary $\partial\tilde N$ such that $\partial{\tilde N}\cap\sigma(A)=\emptyset$ and
$N\cap\sigma(A)\subset\tilde N$, we define the spectral projection $P(A,N)$ by
$$P(A,N):=-\frac{1}{2\pi\sqrt{-1}}\int_{\partial\tilde N}(A-\zeta I)^{-1}d\zeta.$$
The orientation of $\partial\tilde N$ is chosen to make $\tilde N$ stays in the left side of $\partial\tilde N$.

Inspired by \cite{Ph96}, we find that the necessary data for defining the spectral flow are the following:
\begin{itemize}
\item a co-oriented bounded real $1$-dimensional regular $C^1$
submanifold $\ell$ of $\C$ without boundary (we call such an $\ell$ {\bf admissible}, and denote by
$\ell\in\Aa(\C)$);
\item a complex Banach space $X$ (for a real one $X$, we consider $X\bigotimes\C$);
\item and a continuous family (in the gap norm sense) of admissible operators $A_s$, $0\le
s\le 1$ in $\Aa_{\ell}(X)$.
\end{itemize}

Here we define $A\in\Cc(X)$ to be {\bf admissible} with respect to $\ell$, if there exists a bounded open
neighbourhood $N$ of $\ell$ in $\C$ with $C^1$ boundary $\partial N$ such that (i) $\partial
N\cap\sigma(A)=\emptyset$; (ii) $N\cap\sigma(A)\subset \ell$ is a finite set; and (iii) $P^0_{\ell}(A):= P(A,N)$
is a finite rank projection.

We call $\nu_{h,\ell}(A):=\dim \im P^0_{\ell}(A)$ the {\bf hyperbolic nullity} of $A$ with respect to $\ell$. We
denote by $\Aa_{\ell}(X)$ the set of closed admissible operators with respect to $\ell$. It is an open subset of
$\Cc(X)$.

\begin{example}\label{ex05m3.1}
a) In the self-adjoint case, $\ell=\sqrt{-1}(-\epsilon,\epsilon)$ ($\epsilon>0$) with co-orientation from left to
right. Then a self-adjoint operator $A$ is admissible with respect to $\ell$ if and only if $A$ is Fredholm.

b) Another important case is that $\ell=(1-\epsilon,1+\epsilon)$ ($\epsilon\in(0,1)$) with co-orientation from
downward to upward, and all $A_s$ unitary. A unitary operator $A$ is admissible with respect to $\ell$ if and only
if $A-I$ is Fredholm.
\end{example}

Similarly as the definition in \cite{Ph96,ZhLo99}, we can define the spectral flow $\sf_{\ell}\{A_s\}$ as follows.
It counts the number of spectral lines of $A_s$ coming from the negative side of $\ell$ to the non-negative side
of $\ell$.

For each $t\in[0,1]$, there exist bounded open subsets $N_t$, $N_t^{\pm}$ of $\C$ such that
$\sigma(A_t)\cap\partial N_t=\emptyset$, $\sigma(A_t)\cap{\bar\ell}\subset N_t\cap\ell$, $N_t=N_t^+\cup
(N_t\cap\ell)\cup N_t^-$, $N_t^{\pm}$ stays in the positive (negative) side of $\ell$ near $N_t\cap\ell$, and
$P(A_t,N_t)$ is a finite rank projection. Here we denote by $\bar\ell$ the closure of $\ell$ in $\C$. Then
$\sigma(A_t)\cap(\partial N_t\cup(\bar\ell\setminus(N_t\cap\ell)))=\emptyset$. The set $(\partial
N_t\cup(\bar\ell\setminus(N_t\cap\ell))$ is compact since it is a bounded closed set. Since the family $\{A_s\}$,
$0\le s\le 1$ is continuous, there exists a $\delta(t)>0$ for each $t\in[0,1]$ such that
$$\sigma(A_s)\cap(\partial N_t\cup(\bar\ell\setminus(N_t\cap\ell)))
=\emptyset\quad\mbox{for all}\;s\in(t-\delta(t),t+\delta(t))\cap[0,1].$$ Then $\sigma(A_s)\cap{\bar\ell}\subset
N_t\cap\ell$, and
\[
\bigl\{P(A_s,N_t)\bigr\}_{s\in (t-\delta(t),t+\delta(t))\cap[0,1]} \quad \mbox{for fixed}\; t\in[0,1],
\]
is a continuous family of projections. By Lemma I.4.10 in Kato \cite{Ka80}, the operators in the family have the
same rank. Since $[0,1]$ is compact, there exist a partition $0=s_0<\ldots<s_n=1$ and $t_k\in[s_k,s_{k+1}]$,
$k=0,\ldots,n-1$ such that $[s_k,s_{k+1}]\subset(t_k-\delta(t_k),t_k+\delta(t_k))$ for each $k=0,\ldots,n-1$.

\begin{definition}\label{d05m3.1} Let $\ell\in\Aa(\C)$ be admissible
and let $\{A_s\}$, $0\le s\le 1$ be a curve in $\Aa_{\ell}(X)$. The {\bf spectral flow} $\sf_{\ell}\{A_s\}$ of the
family $\{A_s\}$, $0\le s\le 1$ with respect to the curve $\ell$ is defined by
\begin{equation}\label{e05m3.2}\sf_{\ell}\{A_s\}=\sum_{k=0}^{n-1}
\left(\dim \im P(A_{s_k},N^-_{t_k})-\dim \im P(A_{s_{k+1}},N^-_{t_k})\right).
\end{equation}
\end{definition}

The spectral flow has the following properties (cf. \cite{Ph96} and Lemma 2.6 and Proposition 2.2 in
\cite{ZhLo99}).

\begin{proposition}\label{p05m3.1} Let $\ell\in\Aa(\C)$ be admissible
and let $\{A_s\}$, $0\le s\le 1$ be a curve in $\Aa_{\ell}(X)$. Then the spectral flow $\sf_{\ell}\{A_s\}$ is
well-defined, and the following properties hold:
\begin{enumerate}
\item[{\rm (i)}] {\bf Catenation.} Assume $t\in [0,1]$. Then we have
\begin{equation}\label{e05m3.3}
\sf_{\ell}\{A_s;0\le s\le t\}+\sf_{\ell}\{A_s;t\le s\le 1\}=\sf_{\ell}\{A_s;0\le s\le 1\}.
\end{equation}

\item[{\rm (ii)}] {\bf Homotopy invariance.} Let $A(s,t)$,
$(s,t)\in[0,1]\times[0,1]$ be a continuous family in $\Aa_{\ell}(X)$. Then we have
\begin{equation} \label{e05m3.4}
\sf_{\ell}\{A(s,t);(s,t)\in\partial([0,1]\times[0,1])\}=0.\end{equation}

\item[{\rm (iii)}] {\bf Endpoint dependence for Riesz continuity.} Let
$\B^{\sa}(X)$, respectively $\Cc^{\sa}(X)$ denote the spaces of bounded, respectively closed self-adjoint
operators in $X$. Let
\[
\begin{array}{ccccc}
R&:&\Cc^{\sa}&\to&\B^{\sa}(X)\\
\ &\, & A&\mapsto & A(A^2+I)^{-\frac{1}{2}}
\end{array}
\]
denote the {\bf Riesz transformation}. Let $A_s\in\Cc^{\sa}(X)$ for $s\in [0,1]$. Assume that $\{R(A_s)\}$ is a
continuous family. If $m^-(A_0)<+\infty$, then $m^-(A_1)<+\infty$ and we have
\begin{equation}\label{e05m3.5}
\sf\{A_s\}=m^-(A_0)-m^-(A_1).
\end{equation}

\item[{\rm (iv)}] {\bf Product.} Let $\{P_s\}$ be a curve of projections on
$X$ such that $P_sA_s\subset A_sP_s$ for all $s\in[0,1]$. Set $Q_s=I-P_s$. Then we have $P_sA_sP_s\in
\Aa_{\ell}(\im P_s)\subset \Cc(\im P_s)$, $Q_sA_sQ_s\in \Aa_{\ell}(\im Q_s)\subset \Cc(\im Q_s)$, and
\begin{equation} \label{e05m3.6}
\sf_{\ell}\{A_s\}=\sf_{\ell}\{P_sA_sP_s\}+\sf_{\ell}\{Q_sA_sQ_s\}.
\end{equation}

\item[{\rm (v)}] {\bf Bound.} For $A\in\Aa_{\ell}(X)$, there exists a
neighbourhood ${\cal N}$ of $A$ in $\Cc(X)$ such that ${\cal N}\subset\Aa_{\ell}(X)$, and for curves $\{A_s\}$ in
${\cal N}$ with endpoints $A_0=:A$ and $A_1=:B$, the {\bf relative Morse index} $I_{\ell}(A,B):=
-\sf_{\ell}\{A_s,0 \, ; \, \le s\le 1\}$ is well defined and satisfies
\begin{equation}\label{e05m3.7}
0\le I_{\ell}(A,B)\le\nu_{h,\ell}(A)-\nu_{h,\ell}(B).
\end{equation}

\item[{\rm (vi)}] {\bf Reverse orientation.} Let $\hat {\ell}$ denote the
curve $\ell$ with opposite co-orientation. Then we have
\begin{equation}\label{e05m3.8}
\sf_{\ell}\{A_s\}+\sf_{\hat {\ell}}\{A_s\}=\nu_{h,\ell}(A_1)-\nu_{h,\ell}(A_0).
\end{equation}

\item[{\rm (vii)}] {\bf Zero.} Suppose that $\nu_{h,\ell}(A_s)$ is constant for
$s\in [0,1]$. Then $\sf_{\ell}\{A_s\}=0$.

\item[{\rm (viii)}] {\bf Invariance.} Let $\{T_s\}_{s\in[0,1]}$ be a curve of
bounded invertible operators. Then we have
\begin{equation}\label{e05m3.9}\sf_{\ell}\{T_s^{-1}A_sT_s\}=\sf_{\ell}\{A_s\}.
\end{equation}

\end{enumerate}
\end{proposition}

\bp We shall only prove the spectral flow is well-defined. The proof for the rest of the proposition is the same
as that in \cite{Ph96} and Lemma 2.6 and Proposition 2.2 in \cite{ZhLo99} and is omitted.

Since two different partitions of $[0,1]$ has a common refinement, we only need to prove the following local
result:

{\bf Claim.} Let $N_l$, $N_l^{\pm}$, $l=1,2$ be open subsets in $\C$. Assume that for all $s\in[0,1]$ and $l=1,2$,
we have $\sigma(A_s)\cap\partial N_l=\emptyset$, $\sigma(A_s)\cap{\bar\ell}\subset N_l\cap\ell$, $N_l=N_l^+\cup
(N_l\cap\ell)\cup N_l^-$, $N_l^{\pm}$ stays in the positive (negative) side of $\ell$ near $N_l\cap\ell$, and
$P(A_s,N_l)$ is a finite rank projection. Then we have
\[\dim \im P(A_0,N_1^-)-\dim \im P(A_1,N_1^-)
=\dim \im P(A_0,N_2^-)-\dim \im P(A_1,N_2^-).\]

In fact, our assumptions implies
\[\sigma(A_s)\cap\partial(N_1^-\setminus N_2^-)=\sigma(A_s)\cap\partial(N_2^-\setminus N_1^-)=\emptyset.\]
Then $P(A_s,N_1^-\setminus N_2^-)$ and $P(A_s,N_2^-\setminus N_1^-)$, $s\in[0,1]$ are continuous family of
projections. By Lemma I.4.10 in Kato \cite{Ka80}, $\im P(A_t,N_1^-\setminus N_2^-)$ and $\im P(A_t,N_2^-\setminus
N_1^-)$ are constants. So we have
\begin{eqnarray*}& &\left(\dim \im P(A_0,N_1^-)-\dim \im P(A_1,N_1^-)\right)
-\left(\dim \im P(A_0,N_2^-)-\dim \im P(A_1,N_2^-)\right)\\
&=&\left(\dim \im P(A_0,N_1^-)-\dim \im P(A_0,N_2^-)\right)
-\left(\dim \im P(A_1,N_1^-)-\dim \im P(A_1,N_2^-)\right)\\
&=&\left(\dim \im P(A_0,N_1^-\setminus N_2^-)-\dim \im P(A_0,N_2^-\setminus N_1^-)\right)\\
& &-\left(\dim \im P(A_1,N_1^-\setminus N_2^-)-\dim \im P(A_1,N_2^-\setminus N_1^-)\right)\\
&=&0.
\end{eqnarray*}
Thus our claim is proved.
\qe

\begin{remark} \label{r05m3.1} In (iv) of the above proposition, we allow the Banach space $\im P_s$ continuous varying. By \cite[Lemma I.4.10]{Ka80},
for $t\in[0,1]$ being close enough to $s$, there is a continuous family of invertible operators $U_{s,t}\in\B(X)$ such that
\[P_tU_{s,t}=U_{s,t}P_s, \qquad U_{s,t}\to I,\;\mbox{as}\;t\to s.\]
So locally we can define the spectral flow of $B_t\in\Cc(\im P_t)$ as that of $U_{s,t}^{-1}B_tU_{s,t}:\im P_s\to\im P_s$ ($s$ fixed), and globally patch them
together.
\end{remark}

\ss{Calculation of the spectral flow}\label{ss:cal spectral flow}

In this subsection we shall give a method of calculating the spectral flow of differentiable curves, inspired
among others by J.J. Duistermaat \cite{Du76} and J. Robbin and D. Salamon \cite{RoSa93}.

Let $X$ be a complex Banach space, $\tilde N\subset\tilde N$ be bounded open subsets of $\C$, and $\gamma$ be a
closed $C^1$ curve in $\C$ which bounds $\tilde N$. Let $A_s$, $s\in(-\vp,\vp)$, where $\vp>0$, be a curve in
$\C(X)$. Assume that $\gamma\cap\sigma(A_s)=\emptyset$ and $N\cap\sigma(A_s)\subset\tilde N$ for all
$s\in(-\vp,\vp)$. Set $A:=A_0$, $P_s:=P(A_s,N)$, and $P:=P_0$. Assume that $\im P\subset\dom(A_s)$ for all
$s\in(-\vp,\vp)$, $\im P$ is a finitely dimensional subspace of $X$, and $\frac{d}{ds}|_{s=0}(A_sP)=B$ (in the
bounded operator sense). Let $f$ be a polynomial. Then $P_sf(A_s)P_s$, $s\in(-\vp,\vp)$ is a continuous family of
bounded operators, and
\begin{equation}\label{e05m3.10}
P_sf(A_s)P_s=-\frac{1}{2\pi\sqrt{-1}}\int_{\gamma}f(\zeta)(A-\zeta I)^{-1}d\zeta.
\end{equation}
Since $P_s$ $s\in(-\vp,\vp)$ is a continuous family, we have $\|P_s-P\|<1$ if $|s|$ is small. For such $s$, set
$R_s=(I-(P_s-P)^2)^{-\frac{1}{2}}$. Since $P(P_s-P)^2=(P_s-P)^2P$ and $P_s(P_s-P)^2=(P_s-P)^2P_s$, we have
$R_sP=PR_s$ and $R_sP_s=P_sR_s$. Set
\begin{eqnarray*}
U_s^{'}&=&P_sP+(I-P_s)(I-P),\qquad U_s=U_s^{'}R_s,\\
V_s^{'}&=&PP_s+(I-P)(I-P_s),\qquad V_s=V_s^{'}R_s.
\end{eqnarray*}
Then we have
\begin{eqnarray*}
U_sV_s&=&V_sU_s=I,\\
U_sP&=&P_sU_s=P_sR_sP,\\
PV_s&=&V_sP_s=PR_sP_s.
\end{eqnarray*}

\begin{lemma}\label{l05m3.1} We have
\begin{equation}\label{e05m3.11}
\frac{d}{ds}|_{s=0}(U_s^{-1}P_sA_sP_sU_s)=\frac{1}{2\pi\sqrt{-1}}\int_{\gamma}\zeta(A-\zeta I)^{-1}PB(A-\zeta
I)^{-1}d\zeta.
\end{equation}
If $(PAP)(PB)=(PB)(PAP)$, then we have
\begin{eqnarray}\frac{d}{ds}|_{s=0}(P_sP)&=&0,\nn
\\\label{e05m3.12}
\frac{d}{ds}|_{s=0}(U_s^{-1}P_sA_sP_sU_s)&=&PB.
\end{eqnarray}
\end{lemma}

\bp By the definition of $U_s$ and $V_s$ we have
$$U_s^{-1}P_sA_sP_sU_s=V_sP_sA_sP_sU_s=PR_sP_sA_sP_sR_sP.$$
By (\ref{e05m3.11}) we have
\begin{equation}\label{e05m3.13}
(P_sf(A_s)P_s-Pf(A)P)P=\frac{1}{2\pi\sqrt{-1}}\int_{\gamma}f(\zeta)(A_s-\zeta I)^{-1}(A_sP-AP)(A-\zeta
I)^{-1}d\zeta.
\end{equation}
Since $A_s$, $s\in(-\vp,\vp)$ is a curve in $\Cc(X)$ and $\im P$ has finite dimension, we have
\begin{equation}\label{e05m3.14}
\frac{d}{ds}|_{s=0}(P_sf(A_s)P_sP)=\frac{1}{2\pi\sqrt{-1}}\int_{\gamma}f(\zeta)(A-\zeta I)^{-1}B(A-\zeta
I)^{-1}d\zeta.
\end{equation}
Take $f=1$, we have $\frac{d}{ds}|_{s=0}(P_sP)$ exists. By the definition of $R_s$ we have
$\frac{d}{ds}|_{s=0}(R_sP)=0$. Hence we have
\begin{eqnarray*}
\frac{d}{ds}|_{s=0}(U_s^{-1}P_sA_sP_sU_s)
&=&\frac{d}{ds}|_{s=0}(PR_sP_sA_sP_sR_sP)\\
&=&\frac{d}{ds}|_{s=0}((R_sP)(P_sA_sP_sP)(R_sP))\\
&=&\frac{1}{2\pi\sqrt{-1}}\int_{\gamma}\zeta P(A-\zeta I)^{-1}B(A-\zeta I)^{-1}Pd\zeta\\
&=&\frac{1}{2\pi\sqrt{-1}}\int_{\gamma}\zeta(A-\zeta I)^{-1}PB(A-\zeta I)^{-1}d\zeta.
\end{eqnarray*}

In the case of $(PAP)(PB)=(PB)(PAP)$, we have
\[
\frac{d}{ds}|_{s=0}(P_sP)=\frac{1}{2\pi\sqrt{-1}}\int_{\gamma}(A-\zeta I)^{-2}Bd\zeta=0,
\]
and
\begin{eqnarray*} \frac{d}{ds}|_{s=0}(U_s^{-1}P_sA_sP_sU_s)
&=&\frac{1}{2\pi\sqrt{-1}}\int_{\gamma}\zeta P(A-\zeta I)^{-2}Bd\zeta\\
&=&\frac{1}{2\pi\sqrt{-1}}\int_{\gamma}\left(PA(A-\zeta I)^{-2}-P(A-\zeta I)^{-1}\right)Bd\zeta\\
&=&P^2B\\
&=&PB.
\end{eqnarray*}
\qe

\begin{proposition}[cf. Theorem  4.1 of \cite{ZhLo99}]\label{p05m3.2}
Let $X$ be a Banach space and $\ell$ be a bounded open submanifold of $\sqrt{-1}\R$ with co-orientation from left
to right. Let $A_s$, $-\vp\le s\le\vp$ ($\vp>0$), be a curve in $\Aa_{\ell}(X)$. Set $P=P^0_{\ell}(A_0)$, $A=A_0$.
Assume that $\im P\subset\dom(A_s)$ and $B:=\frac{d}{ds}|_{s=0}(A_sP)$ exists. Assume that
\begin{equation}\label{e05m3.15}
(PAP)(PB)=(PB)(PAP),
\end{equation}
where $PAP, PB\in\B(\im P)$, and $PB:\im P\to\im P$ is hyperbolic, i.e. $\sigma(PB)\cap(\sqrt{-1}\R)=\emptyset$.
Then there is a $\delta\in(0,\epsilon)$ such that $\nu_{h,\ell}(A_s)=0$ for all $s\in[-\delta,0)\cup(0,\delta]$
and
\begin{eqnarray}
& &\label{e05m3.16}\sf_{\ell}\{A_s ; 0\le s\le\delta\}=-m^-(PB),\\
& &\label{e05m3.17}\sf_{\ell}\{A_s ; -\delta\le s\le 0\}=m^+(PB).
\end{eqnarray}
Here we denote by $m^+(PB)$ $($$m^-(PB)$$)$ the total algebraic multiplicity of eigenvalues of $PB$ with positive $($negative$)$ imaginary part respectively.
\end{proposition}

\bp We follow the proof of \cite[Theorem 4.1]{ZhLo99}. Since $A\in\Aa_{\ell}(X)$, there exist bounded open subsets
$N$ and $N^{\pm}$ of $\C$ such that $N=N^+\cup(N\cap\ell)\cup N^-$, $N^{\pm}$ stays in the right (left) side of
the imaginary axis, $\sigma(A)\cap{\bar\ell}\subset N\cap\ell$, $\sigma(A)\cap\partial N=\emptyset$, and
$P(A,N)=P$. Since $A_s$, $s\in(-\vp,\vp)$ is a continuous family in $\Cc(X)$, $\sigma(A_s)\cap(\partial N
\cup(\bar\ell\setminus(N\cap\ell)))=\emptyset$ for $|s|$ small. For such $s$, let $P_s$ be defined in Lemma
\ref{l05m3.1}. Then $\|P_s-P\|<1$ for $|s|$ small, and $R_s$ and $U_s$ in Lemma \ref{l05m3.1} are well-defined for
such $s$. Then we have
\[\sigma(A_s)\cap\ell\subset\sigma(A_s)\cap N=\sigma(U_s^{-1}P_sA_sP_sU_s).\]

Now we work in the finite dimensional vector space $\im P$. Since $PB$ commutes with $PAP$, we can assume that
they are both in Jordan normal forms. Then $P(A+sB)P$ is also in Jordan norm form for each $s$. By Lemma
\ref{l05m3.1}, we have $\frac{d}{ds}|_{s=0}(U_s^{-1}P_sA_sP_sU_s)=PB$. Then there exists a $\delta\in(0,\vp)$ such
that $U_s^{-1}P_sA_sP_sU_s$ are hyperbolic for all $s\in[-\delta,0)\cup(0,\delta]$, and
\begin{eqnarray*}
m^-(U_s^{-1}P_sA_sP_sU_s)&=&m^-(PB)\quad\mbox{for all}\;s\in(0,\delta],
\\m^-(U_s^{-1}P_sA_sP_sU_s)&=&m^+(PB)\quad\mbox{for all}\;s\in[-\delta,0).
\end{eqnarray*}
Then our results follows form the definition of the spectral flow and the fact that
\[\dim \im
P(A_s,N^-)=m^-(U_s^{-1}P_sA_sP_sU_s)\quad\mbox{for all}\;s\in[-\delta,0)\cup(0,\delta].\] \qe

\ss{Spectral flow for curves of quadratic forms}\label{ss:quadratic}

Let $X$ be a complex Hilbert space and $\ell=\sqrt{-1}(-\epsilon,\epsilon)$ ($\epsilon>0$) with co-orientation
from left to right. Let $A_s$, $0\le s\le 1$ be a curve of closed self-adjoint Fredholm operators. We will denote
by $\sf\{A_s\}=\sf_{\ell}\{A_s\}$.

\begin{lemma}\label{l05m3.2}
Let $X$ be a Hilbert space. Let $A_s$, $0\le s\le 1$ be a curve of closed self-adjoint Fredholm operators. Then
for any curve $P_s\in\B(X)$ of invertible operators, we have
\begin{equation}\label{e05m3.18}
\sf\{P_sP_s^*A_s\}=\sf\{P_s^*A_sP_s\}=\sf\{A_s\}.
\end{equation}
\end{lemma}

\bp Since $A_s$ is a curve of closed self-adjoint Fredholm operators and $P_s$ is a curve of bounded invertible operators, the families $P_s^*A_sP_s$ and
$P_sP_s^*A_s$, $0\le s\le 1$ are curves of closed Fredholm operators. By (viii) of Proposition \ref{p05m3.1} we have
\begin{eqnarray}
\sf\{P_sP_s^*A_s\}&=&\sf\{P_s(P_s^*A_sP_s)P_s^{-1}\}\nn\\
&=&\sf\{P_s^*A_sP_s\}.\label{e05m3.19}
\end{eqnarray}
Since $P_s^*A_tP_s$ are self-adjoint Fredholm operators and $\dim \ker(P_s^*A_tP_s)=\dim \ker A_t$, we have
\begin{eqnarray}
\sf\{P_s^*A_sP_s\}&=&\sf\{P_0^*A_sP_0\}+\sf\{P_s^*A_1P_s\}\nn\\
&=&\sf\{P_0^*A_sP_0\}\nn\\
&=&\sf\{P_1^*A_sP_1\}.\label{e05m3.20}
\end{eqnarray}

Let $Q_s$ $0\le s\le 1$ be a curve of bounded positive definite operators on $X$ with $Q_0=I$, $Q_1^2=P_0P_0^*$.
By (\ref{e05m3.19}) and (\ref{e05m3.20}) we have
\begin{eqnarray*}
\sf\{P_s^*A_sP_s\}&=&\sf\{P_0^*A_sP_0\}\\
&=&\sf\{P_0P_0^*A_s\}\\
&=&\sf\{Q_1A_sQ_1\}\\
&=&\sf\{Q_0A_sQ_0\}\\
&=&\sf\{A_s\}.
\end{eqnarray*}
\qe

The above lemma leads the following definition.

\begin{definition}\label{d05m3.2} Let $X$ be a Hilbert space. Let $\Ii_s$, $0\le s\le 1$ be a curve of
{\bf bounded Fredholm quadratic forms}, i.e. $\Ii_s(x,y)=\langle A_sx,y\rangle_X$ for all $x,y\in X$, where $A_s$,
$0\le s\le 1$ is a curve of bounded self-adjoint Fredholm operators, and $\langle\cdot,\cdot\rangle_X$ denotes the
inner product in $X$.
\begin{enumerate}
\item[{\rm (a)}] The {\bf spectral flow} $\sf\{\Ii_s\}$ of $\Ii_s$ is defined to be the spectral flow $\sf\{A_s\}$.

\item[{\rm (b)}] If $A_1-A_0$ is compact, the {\bf relative Morse index} $I(\Ii_0,\Ii_1)$ is defined to be the
relative Morse index $I(A_0,A_1):=-\sf\{A_0+s(A_1-A_0)\}$.
\end{enumerate}
\end{definition}

Based on this observation we have the following lemma.

\begin{lemma}\label{l05m3.3} Let $X$ be a Hilbert space. Let $A_s\in\B(X)$, $0\le s\le 1$ be a curve of self-adjoint Fredholm operators and $\Ii_s$ be
quadratic forms defined by $\Ii_s(x,y)=\langle A_sx,y\rangle$ for all $x,y\in X$. Assume that $P_s\in\B(X)$, $0\le s\le 1$ is a curve of operators such that
$P_s^2=P_s$ and $\Ii_s(x,y)=0$ for all $x\in\im P_s$, $y\in\im Q_s$, where $Q_s=I-P_s$. Then we have
\begin{equation}\label{e05m3.21}
\sf\{\Ii_s\}=\sf\{\Ii_s|_{\im P_s}\}+\sf\{\Ii_s|_{\im Q_s}\}.
\end{equation}
\end{lemma}

\bp Set $R_s:=P_s^*P_s+Q_s^*Q_s$, $s\in[0,1]$. Since $P_s+Q_s=I$ and $P_s^2=I$, we have
\[R_s=\frac{I}{2}+2(\frac{I}{2}-P_s^*)(\frac{I}{2}-P_s)>0.\]
Consider the new inner product $\langle R_s x,y\rangle$, $x,y\in X$ on $X$. For this inner product $P_s$ is an orthogonal projection, i.e. $R_sP_s=P_s^*R_s$.

Now we work in the Hilbert space $X$ with the new inner product. So we can assume that $P_s$ is orthogonal. By the fact that $\im P_s$ and $\im Q_s$ are
$\Ii_s$ orthogonal, we have $P_sA_sQ_s=Q_sA_sP_s=0$. Then we have
\[A_s=(P_s+Q_s)A_s(P_s+Q_s)=P_sA_sP_s+Q_sA_sQ_s.\]
So $P_sA_s=A_sP_s$. By (iv) of Proposition \ref{p05m3.1}, $P_sA_sP_s$ is a Fredholm operator on $\im P_s$, $Q_sA_sQ_s$ is a Fredholm operator on $\im Q_s$, and
we have
\begin{eqnarray*}
\sf\{\Ii_s\}&=&\sf\{A_s\}\\
&=&\sf\{P_sA_sP_s:\im P_s\to\im P_s\}+\sf\{Q_sA_sQ_s:\im Q_s\to\im Q_s\}\\
&=&\sf\{\Ii_s|_{\im P_s}\}+\sf\{\Ii_s|_{\im Q_s}\}.
\end{eqnarray*}

\begin{lemma}\label{l05m3.4} Let $X$ be a Hilbert space, and $M$ be a closed subspace with finite codimension. Let $A\in\B(M)$ be a self-adjoint
Fredholm operator and $\Ii(x,y)=\langle Ax,y\rangle$ for all $x,y\in M$. Let $N_1$ and $N_2$ be subspaces of $H$
such that $X=M\oplus N_1=M\oplus N_2$. Define $\Ii_k$ on $H$, $k=1,2$ by
\[\Ii_k(x+u,y+v)=\langle Ax,y\rangle,\qquad\mbox{for all}\; x,y\in M\;\mbox{and}\;u,v\in N_k.\]
Then we have $I(\Ii_1,\Ii_2)=0$.
\end{lemma}

\bp Let $N_0$ be the orthogonal complement of $M$. Set $A_0=\diag(A,0)$ under the direct sum decomposition $X=M\oplus N_0$. Define $\Ii_0$ and $A_1$, $A_2$ by
$\Ii_k(x,y)=\langle A_kx,y\rangle$, for all $x,y\in H$, where $k=0,1,2$. Let $B:N_1\to N$ be a linear isomorphism. Define $P_1\in\B(X)$ by $P_1(x+y)=x+By$ for
all $x\in M$, $y\in N_1$. Then $P_1$ is invertible, $P_1-I$ is compact, and $A_1=P_1^*A_0P_1$. So $A_1-A_0$ is compact. Let $P_s\in\B(X)$, $0\le s\le 1$ be a
curve of invertible operators such that $P_0=I$ and $P_s-I$ are compact. By the definition of the relative Morse index and Lemma \ref{l05m3.2}, we have
\begin{eqnarray*}
I(\Ii_0,\Ii_1)&=&I(A_0,A_1)\\
&=&I(A_0,A_1)\\
&=&-\sf\{P_s^*A_0P_s\}\\
&=&-\sf\{A_0\}\\
&=&0.
\end{eqnarray*}
similarly we have $A_2-A_0$ is compact and $\Ii(A_0,A_2)=0$. So $A_2-A_1$ is compact, and
\[I(\Ii_1,\Ii_2)=I(\Ii_0,\Ii_2)-I(\Ii_0,\Ii_1)=0.\]
\qe

The following proposition gives a generalization of Proposition 5.3 in \cite{AgSa96} and a formula of M. Morse.

\begin{proposition}\label{p05m3.3} Let $X$ be a Hilbert space and $A\in\B(X)$ be a self-adjoint Fredholm operator.
Let $P$ be an orthogonal projection such that $\ker P$ is of finite dimension. Let $\Ii$ be a quadratic form on $X$ defined by $\Ii(x,y)=\langle Ax,y\rangle$,
$x,y\in X$. Set $M=\im P$ and $N$ be the $\Ii$-orthogonal complement of $M$, i.e., $N=\{x\in X;\,\Ii(x,y)=0,\forall y\in M\}$. Then we have
\begin{equation}\label{e05m3.22}
I(PAP,A)=m^-(\Ii|_N)+\dim\ker\Ii|_N-\dim\ker\Ii.
\end{equation}
\end{proposition}

\bp Since $PAP-A$ is of finite rank operator, $sPAP+(1-s)A$, $0\le s\le 1$ is a curve of self-adjoint Fredholm operators. We divide our proof into four steps.

{\bf Step1.} Assume that $\ker A=\{0\}$. Let $M_0=\ker\Ii|_M$, $M_1$ be the orthogonal complement of $M_0$ in $M$,
and $P_0$, $P_1$ be the orthogonal projection onto $M_0$, $M_1$ respectively. Then $P=P_0+P_1$. Since $AM$ is of
finite codimension and $M_0=(AM)^{\perp}\cap M$, $P_0$ is of finite rank. Let $N_1$ be the $\Ii$-orthogonal
complement of $M_1$. Since $M=M_0+M_1$, We have $M_1\cap N_1\subset M_0$. So $M_1\cap N_1=\{0\}$. Moreover we have
\begin{eqnarray*}
\dim N_1&=&\dim\ker(AP_1)-\ind(AP_1)\\
&=&\dim\ker P_1-\ind\, A-\ind\, P_1\\
&=&\dim\ker P_1<+\infty,
\end{eqnarray*}
where we denote $\ind A$ the index of a Fredholm operator $A$. So $X=M_1\oplus N_1$. By the fact that $\Ii$ is nondegenerate, $\Ii|_{N_1}$ is nondegenerate.

Let $\Ii_1$ be defined by $\Ii_1(x+u,y+v)=\Ii(x,y)$ for all $x,y\in M_1$, $u,v\in N_1$. By Lemma \ref{l05m3.3} and Lemma \ref{l05m3.4} we have
\begin{eqnarray*}
I(PAP,A)&=&I(PAP,P_1AP_1)+I(P_1AP_1,A)\\
&=&I(P_1AP_1,A)\\
&=&I(\Ii_1,\Ii)\\
&=&I(\Ii_1|{M_1},\Ii|_{M_1})+I(\Ii_1|{N_1},\Ii|_{N_1})\\
&=&m^-(\Ii|_{N_1}).
\end{eqnarray*}

{\bf Step 2.} Equation (\ref{e05m3.22}) holds if $\ker A=\{0\}$ and $N\subset M$.

In this case, $M_0=N\subset N_1$, $m^-(\Ii|_N)=0$ and $\ker\Ii|N=N$. For each $x\in N_1$ such that $\Ii(x,y)=0$ for all $y\in N$, we have $\Ii(x,y)=0$ for all
$y\in M_1$ and hence for all $y\in M$. Then $x\in N$. Thus $N$ is the $\Ii|_{N_1}$-orthogonal complement of $N$. $N_1$ has an orthogonal decompsition
$N_1=N^+\oplus N^-$ such that $N^+$ and $N^-$ are $\Ii$-orthogonal, $\Ii|_{N^+}>0$ and $\Ii|_{N^-}<0$. Let $P^{\pm}$ be the orthogonal projections onto
$N^{\pm}$. Then $P^{\pm}|_{M_0}$ are isomorphisms. So we have
\[\dim N_1=2\dim N=2m^-(\Ii|_{N_1}).\]
By Step 1 we have
\[I(PAP,A)=m^-(\Ii|_{N_1})=m^-(\Ii|_N)+\dim\ker\Ii|_N-\dim\ker\Ii.\]

{\bf Step 3.} Equation (\ref{e05m3.22}) holds if and $M+N=X$.

In this case we have
\[\ker\Ii|N=\ker\Ii=M\cap N.\]

Firstly we assume that $\ker A=\{0\}$. Then $M_0=\{0\}$, $N_1=N$ and $\ker\Ii|_N=\ker\Ii=\{0\}$. By Step 1, equation (\ref{e05m3.22}) holds.

In the general case, we apply the above special case by taking the quotient space with $\ker A$ and get $I(PAP,A)=m^-(\Ii|N)$.

{\bf Step 4.} Equation (\ref{e05m3.22}) holds.

Firstly we assume that $\ker A=\{0\}$. Let $Q$ be the orthogonal projection onto $M+N$. Then the $\Ii$-orthogonal complement of $M+N$ is $\ker\Ii|_N$. By Step
2 and Step 3 we have
\begin{eqnarray*}
I(PAP,A)&=&I(PAP,QAQ)+I(QAQ,A)\\
&=&m^-(\Ii|N)+\dim\ker\Ii|N.
\end{eqnarray*}

In the general case, we apply the above special case by taking the quotient space with $\ker A$ and get equation (\ref{e05m3.22}). \qe

\ss{A formula}\label{ss:formula}

\begin{lemma}\label{l05m3.5} Let $X$ be a Hilbert space and $H=X\oplus X$. Let $B_s\in\Cc(X)$, $0\le s\le 1$ be a curve of Fredholm
operators. Let the operator $D_s\in\Cc(X)$ by $D_s=\pmatrix{0&B_s^*\cr\\B_s&0\cr}$. Then we have
\begin{equation}\label{e05m3.23}
\sf\{D_s\}=\dim\ker B_1-\dim\ker B_0.
\end{equation}
\end{lemma}

\bp By \cite[Theorem IV.2.23]{Ka80},  $B_s^*$, $0\le s\le 1$ is a curve of closed operators.

Note that $\lambda\in\sigma(D_s)$ if and only if $\lambda^2\in\sigma(B_s^*B_s)$, and the algebraic multiplicities
of them are the same if $|\lambda|\ne 0$ is small. Moreover we have
\begin{eqnarray*}
\dim\ker D_s&=&\dim \ker B_s+\dim\ker B_s^*, \\
\ind B_s=\ind B_0&=&\dim\ker B_s-\dim\ker B_s^*.
\end{eqnarray*}
By the definition of the spectral flow we have
\begin{eqnarray*}
\sf\{D_s\}&=&\frac{1}{2}(\dim\ker D_1-\dim\ker D_0)\\
&=&\dim\ker B_1-\dim\ker B_0.
\end{eqnarray*}
\qe

\begin{lemma}\label{l05m3.6} Let $X$ be a Hilbert space and $H=X\oplus X$. Let $B\in\Cc(X)$ be a operator with compact resolvent,
and $A\in\B(X)$ be a self-adjoint operator. Define linear operator $D_s\in\Cc(X)$ by
$D_s=\pmatrix{sA&B^*\cr\\B&0\cr}$. Then $D_s\in\Cc(H)$, $0\le s\le 1$ is a curve of Fredholm operators, and we
have
\begin{eqnarray}\label{e05m3.24}
\dim\ker D_s&=&\dim\ker A|_{\ker B}+\dim\ker B^*\qquad\mbox{{\rm
for all}}\;s\in(0,1],\\
\label{e05m3.25}\sf\{D_s\}&=&-m^-(A|_{\ker B}).
\end{eqnarray}
\end{lemma}

\bp By \cite[Theorem IV.2.23]{Ka80},  $D_s$, $0\le s\le 1$ is a curve of closed operators. Since $A$ is bounded
and $B$ has compact resolvent, $D_s$ is a Fredholm operator.

For each $s\in(0,1]$ we have
\begin{eqnarray*}
\ker D_s&=&\{(x,y)\in H; sAx+B^*y=0, Bx=0\}\\
&=&\{(x,y)\in H; x\in\ker B, sAx=-B^*y\in\im B^*=\ker B\}\\
&=&\{(x,y)\in H; x\in\ker A|_{\ker B},sAx=-B^*y\}.
\end{eqnarray*}
Define $\varphi:\ker D_s\to\ker A|_{\ker B}$ by $\varphi(x,y)=x$
for $(x,y)\in\ker D_s$. Then $\varphi$ is a linear surjective map,
and $\ker\varphi=\{0\}\times\ker B^*$. Then we get
(\ref{e05m3.24}).

Let $\lambda_t\in\sigma(D_t)$ be a spectral point of $D_t$ near $0$ for $t\ne 0$ small. Then there exists
$(x_t,y_t)\in H\setminus\{0\}$ such that $D_t(x_t,y_t)=\lambda_t(x_t,y_t)$. Then one of the following cases holds.

{\bf Case 1.} $\lambda_t=0$.

In this case, we have $(x_t,y_t)\in\ker D_t$. The algebraic multiplicity of the eigenvalue $0$ of $D_t$ is
$\dim\ker D_t$.

{\bf Case 2.} $\lambda_t\ne 0$ and $Bx_t=0$.

In this case, we have $y_t=0$ and $tAx_t=\lambda_tx_t$. Let $P$ be the orthogonal projection of $X$ onto $\ker B$.
Then $tPAPx_t=\lambda_tx_t$. So the total algebraic multiplicity of these eigenvalues $\lambda_t$ of $D_t$ with
such eigenvectors is
\[m^+(tPAP)+m^-(tPAP)=m^+(A|_{\ker B})+m^-(A|_{\ker B}).\]

{\bf Case 3.} $\lambda_t\ne 0$ and $Bx_t\ne 0$.

In this case, we have $x_t\ne 0$, and
\[\lambda_t^2x_t-t\lambda_tAx_t-B^*Bx_t=0.\]
Take inner product with $x_t$, we have
\begin{equation}\label{e05m3.26}
\lambda_t^2\langle x_t,x_t\rangle-t\lambda_t\langle Ax_t,x_t\rangle-\langle Bx_t,Bx_t\rangle=0.
\end{equation}
For each $x_t$, there exist two $\lambda_t$ satisfying equation
(\ref{e05m3.26}); one is positive, and the other is negative. The
algebraic multiplicity of the two eigenvalues of $D_s$ is equal to
each other. We denote by $2k_t$ the total algebraic multiplicity
of the these eigenvalues of $D_s$ with such eigenvectors.

Since $D_s$, $0\le s\le 1$ is continuous varying, for $t\ne 0$ small, we have
\begin{equation}\label{e05m3.27}
\dim\ker D_0=\dim\ker D_t+m^+(A|_{\ker B})+m^-(A|_{\ker B})+2k_t.
\end{equation}
By the definition of the spectral flow and (\ref{e05m3.24}) we
have
\begin{eqnarray*}
\sf\{D_s\}&=&-m^-(A|_{\ker B})-k_t\\
&=&-m^-(A|_{\ker B})-\frac{1}{2}\left(\dim\ker D_0-\dim\ker D_t-m^+(A|_{\ker B})-m^-(A|_{\ker B})\right)\\
&=&\frac{1}{2}\left(\dim\ker D_1-\dim\ker D_0+\sign(A|_{\ker
B})\right)\\
&=&\frac{1}{2}\left(\dim\ker A|_{\ker B}-\dim\ker B+\sign(A|_{\ker
B})\right)\\
&=&-m^-(A|_{\ker B}).
\end{eqnarray*}
\qe

\begin{proposition}\label{p05m3.4} Let $X$ be a Hilbert space and $H=X\oplus X$. Let $B_s\in\Cc(X)$, $0\le s\le 1$ be a curve of operators with compact resolvent,
and $A_s\in\B(X)$, $0\le s\le 1$ be a curve of self-adjoint operators. Define unbounded operator $D_s$ $X$ by
$D_s=\pmatrix{A_s&B_s^*\cr\\B_s&0\cr}$. Then $D_s\in\Cc(H)$, $0\le s\le 1$ is a curve of Fredholm operators, and we have
\begin{eqnarray}\label{e05m3.28}
\dim\ker D_s&=&\dim\ker A_s|_{\ker B_s}+\dim\ker
B_s^*\qquad\mbox{{\rm
for all}}\;s\in[0,1],\\
\label{e05m3.29}\sf\{D_s\}&=&m^-(A_0|_{\ker B_0})-m^-(A_1|_{\ker
B_1})+\dim\ker B_1-\dim\ker B_0.
\end{eqnarray}
\end{proposition}

\bp (\ref{e05m3.28}) follows form (\ref{e05m3.24}). Set
$D_{s,t}=\pmatrix{tA_s&B_s^*\cr\\B_s&0\cr}$ for $s,t\in[0,1]$. By
\cite[Theorem IV.2.23]{Ka80},  $B_s^*$ and $D_{s,t}$, $0\le s,t\le
1$ are two continuous families of closed operators. Since $A_s$ is
bounded and $B_s$ has compact resolvent, $D_{s,t}$ is a Fredholm
operator.

By Proposition \ref{p05m3.1}, Lemmas \ref{l05m3.5} and \ref{l05m3.6} we have
\begin{eqnarray*}
\sf\{D_s\}&=&-\sf\{D_{0,t};0\le t\le 1\}+\sf\{D_{s,0};0\le s\le 1\}+\sf\{D_{1,t};0\le t\le 1\}\\
&=&m^-(A_0|_{\ker B_0})+(\dim\ker B_1-\dim\ker B_0)
-m^-(A_1|_{\ker B_1})\\
&=&m^-(A_0|_{\ker B_0})-m^-(A_1|_{\ker B_1})+\dim\ker B_1-\dim\ker
B_0.
\end{eqnarray*}
\qe

\s{Maslov-type index theory}\label{maslov-type}

\subsection{Symplectic functional analysis and Maslov index}\label{ss:maslov}

A main feature of symplectic analysis is the study of the {\em Maslov index}. It is an intersection index between
a path of Lagrangian subspaces with the {\em Maslov cycle}, or, more generally, with another path of Lagrangian
subspaces. The Maslov index assigns an integer to each continuous path of Fredholm pairs of Lagrangian subspaces
of a fixed Hilbert space with continuously varying symplectic structures.

Firstly we define symplectic Hilbert spaces and Lagrangian subspaces.

\begin{definition}\label{d05m4.1}
Let $H$ be a complex vector space. A mapping
\[
  \omega:H\times H\to\C
\]
is called a (weak) {\bf symplectic form} on $H$, if it is sesquilinear, skew-symmetric, and non-degenerate, i.e.,
\begin{enumerate}
\item[{\rm (i)}] $\omega(x,y)$ is linear in $x$ and conjugate linear in $y$;

\item[{\rm (ii)}] $\omega(y,x)=-\overline{\omega(y,x)}$;

\item[{\rm (iii)}] $H^{\omega} := \{x\in H \mid \omega(x,y)=0\;\mbox{for all}\; y\in H\} = \{0\}$.
\end{enumerate}
Then we call $(H,\omega)$ a {\bf complex symplectic vector space}.
\end{definition}

\begin{definition}\label{d05m4.2}
Let $(H,\omega)$ be a complex symplectic vector space.
\begin{enumerate}
\item[{\rm (a)}] The {\bf annihilator} of a subspace $\lambda$ of $H$ is defined by
\[
{\lambda}^{\omega} := \{y\in H \mid \omega(x,y)=0\;\mbox{for all}\; x\in\lambda\} .
\]

\item[{\rm (b)}] A subspace $\lambda$ is called {\bf isotropic}, {\bf co-isotropic}, or {\bf Lagrangian} if
\[
\lambda\subset{\lambda}^{\omega}\,,\quad \lambda\supset{\lambda}^{\omega}\,,\quad \lambda={\lambda}^{\omega}
\]
respectively.

\item[{\rm (c)}] The {\bf Lagrangian Grassmannian} $\Ll(H,\omega)$ consists of all Lagrangian subspaces of $(H,\omega)$.
\end{enumerate}
\end{definition}

\begin{definition}\label{d05m4.3}Let $H$ be a complex Hilbert space. A
mapping $\omega:H\times H\to \C$ is called a (strong) {\bf symplectic form} on $H$, if $\omega(x,y)=\langle
Jx,y\rangle_H$ for some bounded invertible skew-symmetric operator $J$. $(H,\omega)$ is called a (strong) {\bf
symplectic Hilbert space}.
\end{definition}

Before giving a rigorous definition of the Maslov index, we fix the terminology and give a simple lemma.

We recall:

\begin{definition}\label{d05m4.4}\begin{enumerate}
\item[{\rm (a)}] The space of (algebraic) {\bf Fredholm pairs} of linear subspaces of a vector space $H$ is defined by
\begin{equation}\label{e05m4.1}
\Ff^{2}_{{\rm alg}}(H):=\{({\lambda},{\mu})\mid \dim\left( {\lambda}\cap{\mu}\right)  <+\infty\; \mbox{and} \dim
\bigl(H/(\lambda+\mu)\bigr)<+\infty\}
\end{equation}
with
\begin{equation}\label{e05m4.2}
\ind(\lambda,\mu):=\dim(\lambda\cap\mu) - \dim(H/(\lambda+\mu)).
\end{equation}

\item[{\rm (b)}] In a Banach space $H$, the space of (topological) {\bf Fredholm pairs} is defined by
\begin{equation}\label{e05m4.3}
\Ff^{2}(H):=\{({\lambda},{\mu})\in\Ff^2_{{\rm alg}}(H)\mid {\lambda},{\mu}, \;\mbox{and}\;{\lambda}+{\mu} \subset
H \;\mbox{is closed}\}.
\end{equation}
\end{enumerate}
\end{definition}

We need the following well-known lemma (see, e.g., \cite[Lemma 1.7]{BoZh04}).

\begin{lemma}\label{l05m4.1}
Let $(H,\omega)$ be a (strong) symplectic Hilbert space. Then
\begin{enumerate}
\item[{\rm (i)}] there exists a $\omega$-orthogonal splitting
\[H=H^+\oplus H^-\]
such that $-\sqrt{-1}\omega$ is positive (negative) definite on $H^{\pm}$, and we call it a {\bf symplectic splitting};

\item[{\rm (ii)}] there is a $1$-$1$ correspondence between the space
$$\Uu(H^+,H^-,\omega)=\{U\in\B\bigl(H^+,H^-\bigr)|\,\,\omega(Ux,Uy)=-\omega(x,y),\forall x,y\in H^+\}$$
and $\Ll(H,\omega)$ under the mapping $U\to L:= \GGG(U)$ $($= graph of $U$$)$;

\item[{\rm (iii)}] if $U,V\in\Uu(H^+,H^-,\omega)$ and $\lambda:=\GGG(U)$, $\mu:=\GGG(V)$, then $(\lambda,\mu)$ is a Fredholm pair if and only if $U-V$, or,
equivalently, $UV^{-1}-I$ is Fredholm. Moreover, we have a natural isomorphism
\begin{equation}\label{e05m4.4}
\ker(UV^{-1}-I) \simeq \lambda\cap \mu\,.
\end{equation}
\end{enumerate}

\end{lemma}

\begin{definition}\label{d05m4.5} Let $(H,\langle\cdot,\cdot\rangle_s)$,
$s\in[0,1]$ be a continuous family of Hilbert spaces, and $\omega_s(x,y)=\langle J_sx,y\rangle_s$ be a continuous
family of symplectic forms on $H$, i.e., $\{A_{s,0}\}$ and $\{J_s\}$ are two continuous families of bounded
invertible operators, where $A_{s,0}$ is defined by
$$\langle x,y\rangle_s=\langle A_{s,0}x,y\rangle_0\quad\mbox{for all}\;x,y\in H.$$
Let $\{(\lambda_s,\mu_s)\}$ be a continuous family of Fredholm pairs of Lagrangian subspaces of $(H,\langle\cdot,\cdot\rangle_s,\omega_s)$. Then there is a
continuous families of symplectic splitting
\begin{equation}\label{e05m4.5}H=H_s^+\oplus
H_s^-\end{equation} for all $s\in[0,1]$. Such $H_s^{\pm}$ can be chosen to be the positive (negative) space associated to the self-adjoint operator
$\sqrt{-1}J_s\in\B(H,\langle\cdot,\cdot\rangle_s)$. By Lemma \ref{l05m4.1}, $\lambda_s=\GGG_s(U_s)$ and $\mu_s=\GGG_s(V_s)$ with $U_s$, $V_s\in
\Uu(H_s^+,H_s^-,\omega_s)$, where $\GGG_s$ denotes the graph associated to the splitting (\ref{e05m4.5}). We define the {\bf Maslov index}
$\Mas\{\lambda_s,\mu_s\}$ by
\begin{equation}\label{e05m4.6}\Mas\{\lambda_s,\mu_s\}=-\sf_{\ell}\{U_sV_s^{^{-1}}\},
\end{equation}
where $\ell:=(1-\epsilon,1+\epsilon)$ with $\epsilon\in(0,1)$ and with upward co-orientation.
\end{definition}

\begin{remark}\label{r05m4.1}
For finite-dimensional $(H,\omega)$, constant $\mu_s=\mu_0$, and a loop $\{\lambda_s\}$, i.e., for $\lambda_0=\lambda_1$\,, we notice that
$\Mas\{\lambda_s,\mu_s\}$ is the winding number of the closed curve $\{\det(U_s^{-1} V_0)\}_{s\in [0,1]}$\,. This is the original definition of the Maslov
index as explained in Arnol'd, \cite{Ar67}.
\end{remark}

\begin{lemma}\label{l05m4.2}The Maslov index is
independent of the choice of the symplectic splitting of $H$.
\end{lemma}

\bp Let $H=H_{s,k}^+\oplus H_{s,k}^-$, $s\in[0,1]$ with $k=0,1$ be two continuous families of symplectic splitting. For each $s\in[0,1]$ and $k=0,1$, set
\[\langle\cdot,\cdot\rangle_{s,k}=(-\sqrt{-1}\omega|_{H_{s,l}^+})\oplus (\sqrt{-1}\omega|_{H_{s,l}^-}),\]
Then $(H,\langle\cdot,\cdot\rangle_{s,k})$ is a Hilbert space for each $s\in[0,1]$ and $k=0,1$. Set
\[\langle\cdot,\cdot\rangle_{s,t}=(1-t)\langle\cdot,\cdot\rangle_{s,0}+t\langle\cdot,\cdot\rangle_{s,1}\]
for each $(s,t)\in[0,1]\times[0,1]$. for each $(s,t)\in[0,1]\times[0,1]$, define $J_{s,t}\in\B(H)$ by
\[\omega(x,y)_s=\langle Jx,y\rangle_{s,t} \quad\mbox{for all}\;x,y\in H.\]
Then $H_{s,k}^{\pm}$ is the positive (negative) space associated the self-adjoint operator $-\sqrt{-1}J_{s,k}$ for each $s\in[0,1]$ and $k=0,1$. Let
$H_{s,t}^{\pm}$ be the positive (negative) space associated the self-adjoint operator $-\sqrt{-1}J_{s,t}$ for each $s\in[0,1]$ and $t=0,1$.

Let $(\lambda_s,\mu_s)$ be a continuous family of Fredholm pairs of Lagrangian subspaces of $(H,\omega_s)$. For each symplectic splitting $H=H_{s,t}^+\oplus
H_{s,t}^-$, we denote by $U_{s,t}$ and $V_{s,t}$ the associated generated "unitary" operators of $\lambda_s$ and $\mu_s$ respectively. We also denote by
$\Mas_t$ the Maslov index defined with $\langle\cdot,\cdot\rangle_{s,t}$ for each $t\in[0,1]$. By Proposition \ref{p05m3.1} we have
\begin{eqnarray*}
\Mas_0\{\lambda_s,\mu_s\}&-&\Mas_1\{\lambda_s,\mu_s\}\\
&=&-\sf_{\ell}\{U_{s,0}V_{s,0}^{^{-1}}\}+\sf_{\ell}\{U_{s,1}V_{s,1}^{^{-1}}\}\\
&=&-\sf_{\ell}\{U_{s,t}V_{s,t}^{^{-1}};(s,t)\in\partial\bigl([0,1]\times[0,1]\bigr)\}\\
&=&0.
\end{eqnarray*}
\qe

\begin{corollary}[Symplectic invariance]\lb{c05m4.1} Let $(H_l,\omega_{s,k})$, $k=1,2$ be two continuous families of symplectic Hilbert spaces.
Let $M(s)\in\B(H_1,H_2)$, $0\le s\le 1$ be a curve of invertible operators such that
\[\omega_{s,2}(M_sx,M_sy)=\omega_{s,1}(x,y)\quad\mbox{for all}\;x,y\in H_1\;\mbox{and}s\in[0,1].\]
Then for any curve $(\lambda(s),\mu(s))$, $0\le s\le 1$ curve of Fredholm pairs of Lagrangian subspaces of $H_1$,
\begin{equation} \label{e05m4.7} \Mas\{M\lambda,M\mu\}=\Mas\{\lambda,\mu\}.
\end{equation}
\end{corollary}

\bp Let $H_1=H^+_{s,1}\oplus H^-_{s,2}$ be a continuous family of symplectic splitting of the family $(H_1,\omega_{s,1})$, $0\le s\le 1$. Then
$H_2=H^+_{s,2}\oplus H^-_{s,2})$ be a continuous family of symplectic splitting of the family $(H_2,\omega_{s,2})$, $0\le s\le 1$, where
$H^+_{s,2}=M_sH^+_{s,1})$ and $H^-_{s,2}=M_sH^-_{s,1}$. For each symplectic splitting $H_k=H^+_{s,k}\oplus H^-_{s,k}$, $s\in[0,1]$ and $k=1,2$, we denote by
$U_{s,k}$ and $V_{s,k}$ the associated generated "unitary" operators of $\lambda_s$ and $\mu_s$ respectively. Then we have
\[U_{s,2}=M_sU_{s,1}M_s^{-1},\qquad V_{s,2}=M_sV_{s,1}M_s^{-1}.\]
By the definition of the Maslov index we have
\begin{eqnarray*} \Mas\{M\lambda,M\mu\}&=&-\sf_{\ell}\{(M_sU_{s,1}M_s^{-1})(M_sV_{s,1}M_s^{-1})^{-1};0\le s\le 1\}\\
&=&-\sf_{\ell}\{U_{s,1}V_{s,1}^{-1};0\le s\le 1\}\\
&=&\Mas\{\lambda,\mu\}.\end{eqnarray*}\qe

Now we give a method of using the crossing form to calculate Maslov indices (cf. \cite{Du76}, \cite{RoSa93} and \cite[Theorem 2.1]{BoFu98}.

Let $\lambda=\{\lambda_s\}_{s\in[0,1]}$ be a $C^1$ curve of Lagrangian subspaces of $(H,\omega)$. Let $W$ be a fixed Lagrangian complement of $\lambda_t$. For
$v\in\lambda_t$ and $|s-t|$ small, define $w(s)\in W$ by $v+w(s)\in\lambda_s$. The form
$$
Q(\lambda,t):= Q(\lambda,W,t)(u,v)=\frac{d}{ds}|_{s=t}\omega(u,w(s)),\quad \forall u,v\in\lambda_t
$$
is independent of the choice of $W$. Let $\{(\lambda_s,\mu_s)\}$, $0\le s\le 1$ be a curve of Fredholm pairs of
Lagrangian subspaces of $H$. For $t\in[0,1]$, the {\bf crossing form} $\Gamma(\lambda,\mu,t)$ is a quadratic form
on $\lambda_t\cap\mu_t$ defined by
$$
\Gamma(\lambda,\mu,t)(u,v)=Q(\lambda,t)(u,v)-Q(\mu,t)(u,v),\quad\forall u,v\in\lambda_t\cap\mu_t\/.
$$
A {\bf crossing} is a time $t\in[0,1]$ such that $\lambda_t\cap\mu_t\ne\{0\}$. A crossing is called {\bf regular}
if $\Gamma(\lambda,\mu,t)$ is nondegenerate. It is called {\bf simple} if it is regular and $\lambda_t\cap\mu_t$
is one-dimensional.

Now let $(H,\omega)$ be a symplectic Hilbert space with $\omega(x,y)=\langle jx,y\rangle$, for all $x,y\in H$,
where $J\in\B(H)$ is a invertible skew self-adjoint operator. Then we have a symplectic Hilbert space $X=(H\oplus
H,(-\omega)\oplus\omega)$. For each $M\in\Sp(H,\omega)$, its graph $\Gr(M)$ is a Lagrangian subspace of $X$. The
following lemma is Lemma 3.1 in \cite{Du76}.

\begin{lemma} \label{l05m4.3} Let $M(s)\in\Sp(H,\omega)$, $0\le s\le 1$ be a curve of linear symplectic maps. Assume that
$M(s)$ is differentiable at $t\in[a,b]$. Set $B_1(t)=-J{\dot M(t)}M(t)^{-1}$ and $B_2(t)=-JM(t)^{-1}{\dot M(t)}$.
Then $B_1(t)$, $B_2(t)$ are self-adjoint, $B_2(t)=M(t)^*B_1(t)M(t)$ and we have

\begin{equation} \label{e05m4.8}
Q(\Gr(M),t)((x,M(t)x),(y,M(t)y))=(B_2(t)x,y).
\end{equation}\end{lemma} \qe



\begin{proposition} \label{p05m4.1} Let $(H,\omega)$ be a symplectic Hilbert space
and $\{(\lambda_s,\mu_s)\}$, $0\le s\le 1$ be a $C^1$ curve of Fredholm pairs of Lagrangian subspaces of $H$ with
only regular crossings. Then we have
\begin{equation} \label{e05m4.9}
\Mas\{\lambda,\mu\}=m^+(\Gamma(\lambda,\mu,0))
-m^-(\Gamma(\lambda,\mu,1))+\sum_{0<t<1}\sign(\Gamma(\lambda,\mu,t)).
\end{equation}
\end{proposition}

\bp Pick an invertible skew self-adjoint operator $J\in\B(H)$ such that $J^2=-I$ and $\omega(x,y)=\langle
Jx,y\rangle$. Let $H_1=\ker(J-\sqrt{-1}I)$ and $H_2=\ker(J+\sqrt{-1}I)$. By Lemma \ref{l05m4.1} there are curves
of isometric $U(t)$, $V(t)$ in $\U(H_1,H_2)$ such that $\lambda(t)=\Gr(U(t))$ and $\mu(t)=\Gr(V(t))$. Apply Lemma
\ref{l05m4.3} for $(H_1,\langle -\sqrt{-1}x,y\rangle)$, for any $x,y\in\ker(U(t)-V(t))$ and $t\in[a,b]$ we have
\begin{eqnarray*}
\frac{d}{ds}|_{s=t}\langle-\sqrt{-1}V^{-1}Ux,y\rangle
&=&\langle\sqrt{-1}V^{-1}{\dot V}V^{-1}Ux,y\rangle+\langle-\sqrt{-1}V^{-1}{\dot U}x,y\rangle\\
&=&\langle\sqrt{-1}V^{-1}{\dot V}x,y\rangle+\langle-\sqrt{-1}U^{-1}VV^{-1}{\dot U}x,U^{-1}Vy\rangle\\
&=&\langle\sqrt{-1}V^{-1}{\dot V}x,y)+\langle-\sqrt{-1}U^{-1}{\dot U}x,y\rangle\\
&=&-\Gamma(\lambda,\mu,t)((x,Ux),(y,Uy)).
\end{eqnarray*}

By Proposition \ref{p05m3.2} we obtain (\ref{e05m4.7}). \qe

\ss{Spectral flow formula for fixed maximal domain}\label{ss:spectral-fixd}

Let $D_m\hookrightarrow D_M\hookrightarrow X$ be three Hilbert spaces. We assume that $D_m$ is a closed subspace
of $D_{M}$ and a dense subspace of $X$. Let $\{A_s\}_{s\in[0,1]}$ be a family of symmetric densely defined
operators in $\Cc(X)$ with domain $\dom(A_s)=D_m$. Here we denote by $\Cc(X)$ all closed operators in $X$. Assume
that $\dom(A_s^*)=D_{M}$\,, i.e., the domain of the maximal symmtric extension $A^*_s$ of $A_s$ is independent of
$s$.

We recall from \cite{BoFu98} for each $s\in[0,1]$:

\begin{enumerate}

\item The space $D_{M}$ is a Hilbert space with the graph inner
product
\begin{equation}\label{e05m4.10}
\langle x,y\rangle_{\GGG_s} := \langle x,y\rangle_X + \langle A_s^*x,A_s^*y\rangle_X \quad\mbox{for}\; x,y\in
D_M\,.
\end{equation}

\item The space $D_m$ is a closed subspace in the graph norm and
the quotient space $D_M/D_m$ is a strong symplectic Hilbert space with the (bounded) symplectic form induced by
Green's form
\begin{equation}\label{e05m4.11}
\omega_s(x+D_m,y+D_m) := \langle A_s^*x,y\rangle_X - \langle x,A_s^*y\rangle_X \quad\mbox{for}\; x,y\in D_M\,.
\end{equation}

\item If $A_s$ admits a self-adjoint Fredholm extension
$A_{s,D_s}:=A_s^*|_{D_s}$ with domain $D_s\subset X$, then the {\bf natural Cauchy data space} $(\ker
A_s^*+D_m)/D_m$ is a Lagrangian subspace of $(D_M/D_m,\omega_s)$\,.

\item Moreover, self-adjoint Fredholm extensions are characterized
by the property of the domain $D_s$ that $(D_s+D_m)/D_m$ is a Lagrangian subspace of $(D_M/D_m,\omega_s)$ and
forms a Fredholm pair with $(\ker A_s^*+D_m)/D_m$\,.

\item We denote the natural projection (which is independent of
$s$) by
\[
\gamma:D_M\rightarrow D_M/D_m\/.
\]
We call $\gamma$ the {\bf abstract trace map}.
\end{enumerate}

We have the following spectral flow formula (cf. \cite[Theorem 5.1]{BoFu98}, \cite[Corollary 2.14]{BoZh04} and
\cite[Theorem 1.5]{BoZh05}).

\begin{proposition}\label{p05m4.2}
We assume that on $D_M$ the graph norms induced by $A_s^*$ and the original norm are equivalent.
Assume that $\{A_s^*:D_M\to X\}$ is a continuous family of bounded operators and each $A_s$ is injective. Let
$\{D_s/D_m\}$ be a continuous family of Lagrangian subspaces of $(D_M/D_m,\omega_s)$, such that each $A_{s,D_s}$
is a Fredholm operator. Then:

\begin{enumerate}

\item[{\rm (a)}] Each $\bigl(D_s/D_m,\gamma(\ker(A_s^*))\bigr)$ is a Fredholm pair in $D_M/D_m$\/.

\item[{\rm (b)}] Each Cauchy data space $\gamma(\ker A_s^*)$ is a Lagrangian subspace of $(D_M/D_m,\omega_s)$\,.

\item[{\rm (c)}] The family $\{\gamma(\ker A_s^*)\}$ is a continuous family in $D_M/D_m$\,.

\item[{\rm (d)}] The family $\bigl\{A_{s,D_s}\bigr\}$ is a continuous family of self-adjoint Fredholm operators in
$\Cc(X)$.

\item[{\rm (e)}] Finally, we have
\begin{equation}\label{e05m4.12}
\sf\{A_{s,D_s}\} = -\Mas \{\gamma(D_s),\gamma(\ker A_s^*)\}.
\end{equation}
\end{enumerate}
\end{proposition}\qe

\ss{The Maslov-type indices}\label{ss:maslov-type}

\begin{definition} \label{d05m4.6} Let $(X_l,\omega_l)$ be symplectic Hilbert spaces with $\omega_l(x,y)=(j_lx,y)$, $x,y\in X_l$,
$j_l\in\B(X)$ are invertible, and $j_l^*=-j_l$, where $l=1,2$. Then we have a symplectic Hilbert space
$(H=X_1\oplus X_2,(-\omega_1)\oplus \omega_2)$. Let $W\in{\cal L}(H)$. Let $M(t)$, $a\le t\le b$ be a curve in
$\Sp(X_1,X_2)$ such that $\Gr(M(t))\in{\cal FL}(W)$ for all $t\in[a,b]$. The {\bf Maslov-type index} $i_W(M(t))$
is defined to be $\Mas(\Gr(M(t), W)$. If $a=0$, $b=T$, $(X_1,\omega_1)=(X_2,\omega_2)$ and $M(0)=I$, we denote by
$\nu_{T,W}(M(t))=\dim (\Gr(M(T)\cap W)$. \end{definition}

The Maslov-type indices have the following property.

\begin{lemma} \label{l05m4.4}
Let $(X_l,\omega_l)$ be symplectic Hilbert spaces with $\omega_l(x,y)=(j_lx,y)$, where $x,y\in X_l$,
$j_l\in\B(X_l)$ are invertible, and $j_l^*=-j_l$, $l=1,2,3,4$. Let $W$ be a Lagrangian subspace of $(X_1\oplus
X_4,(-\omega_1)\oplus\omega_4)$. Let $\gamma_l\in C([0,1],\Sp(X_l,X_{l+1}))$, $l=1,2,3$ be syplectic paths such
that $\Gr(\gamma_3(s)\gamma_2(t)\gamma_1(s))\in{\cal FL}(W)$ for all $(s,t)\in[0,1]\times[0,1]$. Then we have

\begin{equation}\label{e05m4.13}
i_W(\gamma_3\gamma_2\gamma_1)=i_{W^{'}}(\gamma_2)+i_W(\gamma_3\gamma_2(0)\gamma_1), \end{equation}
where
$W^{'}=\diag(\gamma_1(1),\gamma_3(1)^{-1})W$. \end{lemma}

\bp Let $M=\diag(\gamma_1(1),\gamma_3(1)^{-1})$. By the homotopic
invariance rel. endpoints of the Maslov-type indices and Corollary
\ref{c05m4.1}, we have

\begin{eqnarray*}
i_W(\gamma_3\gamma_2\gamma_1)
&=&i_W(\gamma_3(1)\gamma_2\gamma_1(1))+i_W(\gamma_3\gamma_2(0)\gamma_1)\\
&=&\Mas(M\Gr(\gamma_3(1)\gamma_2\gamma_1(1)),M W)+i_W(\gamma_3\gamma_2(0)\gamma_1)\\
&=&i_{W^{'}}(\gamma_2)+i_W(\gamma_3\gamma_2(0)\gamma_1).
\end{eqnarray*}\qe

The following properties of fundamental solutions for linear ODE
will be used later.

\begin{lemma} \label{l05m4.5} Let $j\in C^1([0,+\infty),\GL(m,\C))$ be a curve of skew self-adjoint matrices, and $b\in
C([0,+\infty),\gl(m,\C))$ be a curve of self-adjoint matrices. Let $\gamma\in C^1([0,+\infty),\GL(m,\C))$ be the
fundamental solution of

\begin{equation} \label{e05m4.14} -j{\dot x}-\frac{1}{2}{\dot j}x=bx. \end{equation}
Then we have $\gamma(t)^*j(t)\gamma(t)=j(0)$ for all $t$. \end{lemma}

\bp By the definition of the fundamental solution, we have
$\gamma(0)^*j(0)\gamma(0)=j(0)$. Since $j^*=-j$ and $b^*=b$, we have

\begin{eqnarray*}
\frac{d}{dt}(\gamma(t)^*j(t)\gamma(t)) &=&{\dot \gamma }^*j\gamma+\gamma^*{\dot
j}\gamma+\gamma^*j{\dot\gamma}\\
&=&(-b\gamma-\frac{1}{2}{\dot j})^*j^{*-1}j\gamma+\gamma^*{\dot
j}\gamma+\gamma^*jj^{-1}(-b\gamma-\frac{1}{2}{\dot j})\\
&=&\gamma^*(b-\frac{1}{2}{\dot j}+{\dot j}-b-\frac{1}{2}{\dot
j})\gamma\\
&=&0.
\end{eqnarray*}
So we have $\gamma(t)^*j(t)\gamma(t)=j(0)$.
\qe

\begin{lemma} \label{l05m4.6} Let $B\in C([0,+\infty),\gl(m,\C))$ and $P\in C^1([0,+\infty),\GL(m,\C))$ be two curves of
matrices. Let $\gamma\in C^1([0,+\infty),\GL(m,\C))$ be the fundamental solution of
\begin{equation} \label{e05m4.15} {\dot x}=Bx, \end{equation}
and $\gamma^{'}\in C^1([0,+\infty),\GL(m,\C))$ be the fundamental solution of
\begin{equation} \label{e05m4.16}
{\dot y}=(PBP^{-1}+{\dot P}P^{-1})y.
\end{equation}
Then we have
\begin{equation} \label{e05m4.17}
\gamma^{'}=P\gamma P(0)^{-1}.
\end{equation}
\end{lemma}

\bp Direct calculation shows

$$\frac{d}{dt}(P\gamma P(0)^{-1})=(PBP^{-1}+{\dot P}P^{-1})P\gamma
P(0)^{-1}$$ and $P(0)\gamma P(0)^{-1}=I$. By definition, $P\gamma P(0)^{-1}$ is the fundamental solution of
(\ref{e05m4.16}). \qe

\begin{corollary} \label{c05m4.2}
Let $j_1,j_2\in C^1([0,+\infty),\GL(m,\C))$ be two curves of skew self-adjoint matrices. Let
$P\in C^1([0,+\infty),\GL(m,\C))$ be a curve of matrices such that $P^*j_2P=j_1$, and $b\in
C([0,+\infty),\GL(m,\C))$ be a curve of self-adjoint matrices. Let $\gamma\in C^1([0,+\infty),\GL(m,\C))$ be the
fundamental solution of
\begin{equation} \label{e05m4.18}
-j_1{\dot x}-\frac{1}{2}{\dot  j}_1x=bx, \end{equation}
and $\gamma^{'}\in C^1([0,+\infty),\GL(m,\C))$ be the
fundamental solution of
\begin{equation} \label{e05m4.19} j_2{\dot y}-\frac{1}{2}{\dot
j}_2y=(P^{*-1}bP^{-1}+Q)y, \end{equation}
where $Q=\frac{1}{2}(P^{*-1}{\dot P}^*j_2-j_2{\dot P}P^{-1})$. Then we
have

\begin{equation} \label{e05m4.20}
\gamma^{'}=P\gamma P(0)^{-1}.
\end{equation}

In particular, when $j_1$ and $j_2$ are constant matrices, we have
$$Q=P^{*-1}{\dot P}^*j_2=-j_2{\dot P}P^{-1}.$$
\end{corollary}

\bp Take $B=-j_1^{-1}(b+\frac{1}{2}{\dot  j}_1)$ in Lemma \ref{l05m4.6}, we have
\begin{eqnarray*}
-j_2(PBP^{-1}+{\dot P}P^{-1})-\frac{1}{2}{\dot  j}_2 &=&-j_2(P(-j_1)^{-1}(b+\frac{1}{2}{\dot j}_1)P^{-1}+{\dot
P}P^{-1})-\frac{1}{2}{\dot
j_2}\\
&=&P^{*-1}(b+\frac{1}{2}{\dot  j}_1)P^{-1}-j_2{\dot
P}P^{-1}-\frac{1}{2}{\dot  j}_2\\
&=&P^{*-1}bP^{-1}-j_2{\dot P}P^{-1}+\frac{1}{2}(P^{*-1}{\dot
j}_1P^{-1}-j_2)\\
&=&P^{*-1}bP^{-1}-j_2{\dot
P}P^{-1}+\frac{1}{2}(P^{*-1}\frac{d}{dt}(P^*j_2P)
P^{-1}-j_2)\\
&=&P^{*-1}bP^{-1}+Q.
\end{eqnarray*}
By Lemma \ref{l05m4.6}, our results holds.\qe

The following is a special case of the spectral flow formula.

Let $j\in C^1([0,T],\GL(m,\C))$ be a curve of skew self-adjoint matrices. Then we have symplectic Hilbert spaces
$(\C^m,\omega(t))$ with standard Hermitian inner product and $\omega(t)(x,y)=(j(t)x,y)$, for all $x,y\in \C^m$ and
$t\in[0,T]$. Then we have a symplectic Hilbert space $(V=\C^m\oplus \C^m,(-\omega(0))\oplus \omega(T))$. Let
$W\in{\cal L}(V)$. Let $b_s(t)\in\B(\C^m)$, $0\le s \le 1$, $0\le t\le T$ be a continuous family of self-adjoint
matrices such that $b_0(t)=0$. By Lemma \ref{l05m4.5}, there are continuous family of matrices
$M_s(t)\in\GL(m,\C)$ such that $M_s(0)=I$, $M_s(t)^*j(t)M_s(t)=j(0)$ and
$$-j\frac{d}{d t}M_s(t)-\frac{1}{2}(\frac{d}{d t}j)M_s(t)=b_s(t)M_s(t).$$
Set
\begin{eqnarray*}
X&=&L^2([0,T],\C^m), \quad D_m=H^1_0([0,T],\C^m), \\
D_M&=&H^1([0,T],\C^m),\quad D_W=\{x\in D_M;(x(0),x(t))\in W\}.
\end{eqnarray*}
Let $A_M\in{\cal C}(X)$ with domain $D_M$ be defined by
$$A_Mx=-j\frac{d}{d t}x-\frac{1}{2}(\frac{d}{d t}j)x.$$
Set $x\in D_M$, $A=A_M|_{D_m}$, $A_W=A_M|_{D_W}$. Let $C_s\in\B(X)$ be defined
by $(C_sx)(t)=b_s(t)x(t)$, $x\in X$, $t\in[0,T]$.

\begin{proposition} \label{p05m4.3} Set $W^{'}=\diag(I,M_0(T)^{-1})W$. Then we have

\begin{equation} \label{e05m4.21}
I(A_W,A_W-C_1)=i_{W^{'}}(M_0^{-1}M_1).
\end{equation}
\end{proposition}

\bp The Sobolev embedding theorem shows that $D_M\subset C([0,T],\C^m)$. For any $x\in D_M$, define
$\gamma(x)=(x(0),x(T))$. Direct calculation shows that $D_M/D_m=\C^m\oplus\C^m$ with symplectic structure
$(\diag(j(0),-j(T))\gamma(x),\gamma(y))$, $x,y\in D_M$, and $\gamma$ is the abstract trace map. Moreover,
$A^*=A_M$, $\gamma(A^*-C_s)=\Gr(M_s(T))$, and $\gamma(D_W)=W$. By Proposition \ref{p05m4.2} and Lemma
\ref{l05m4.4}, we have
\begin{eqnarray*}
I(A_W,A_W-C_1)
&=&-\sf\{A_W-C_s\}\\
&=&\Mas(\{\Gr(M_s(T));0\le s\le 1\}, W)\\
&=&i_W(M_0(T)(M_0(T)^{-1}M_s(T))I;0\le s\le 1)\\
&=&i_{W^{'}}(M_0(T)^{-1}M_s(T);0\le s\le 1)\\
&=&-i_{W^{'}}(M_0(t)^{-1}M_0(t);0\le t\le T)+i_{W^{'}}(M_0(0)^{-1}M_s(0);0\le s\le 1)\\
& &+i_{W^{'}}(M_0(t)^{-1}M_1(t);0\le t\le T)\\
&=&i_{W^{'}}(M_0^{-1}M_1). \end{eqnarray*} \qe

\s{Proof of the main results}\label{s:proof}

In this section we will use the notations in {\S}\ref{s:main}.

\ss{Proof of Theorem \ref{t05m2.1}}\label{ss:pf-t2.1}
\begin{lemma} \label{l05m5.1}
The index forms $\Ii_{s,R}$, $0\le s\le 1$ is a curve of bounded Fredholm quadratic forms on $H_R$.
\end{lemma}

\bp Since $\Ii_{s,R}$ are bounded symmetric quadratic forms on $H_R$, by Riesz representation theorem, they form a
continuous curve.

For each $k,l=0,\ldots,m$ and $s\in[0,1]$, we define the bounded operators $P_{k,l}(s)\in\B(H_R)$ by
\[\langle P_{k,l}(s)x,y\rangle_m=\int_0^T\langle p_{k,l}(s,t)\frac{d^k x}{dt^k},\frac{d^l y}{dt^l}\rangle
dt\quad\mbox{for all}\;x,y\in H_R.\]

{\bf Claim.} $P_{k,l}(s)$ is compact for either $k\ne m$ or $l\ne m$.

Since $P_{k,l}(s)=P_{l,k}^*(s)$, without loss of generality we can assume that $k\ne m$. Pick a bounded sequence
$\{x_{\alpha};\alpha\in\N\}$ in $H_R$. By Sobolev embedding theorem, the sequence $\{p_{k,l}(s,t)\frac{d^k
x_{\alpha}}{dt^k}\}$ has a convergent subsequence, which is denoted by the original sequence. Since $P_{k,l}(s)$
is bounded, we have
\[\lim_{\alpha,\beta\to+\infty}\| P_{k,l}(s)(x_{\alpha}-x_{\beta})\|_m^2
=\lim_{\alpha,\beta\to+\infty}\int_0^T\langle p_{k,l}(s,t)\frac{d^k(x_{\alpha}-x_{\beta})}{dt^k},\frac{d^l
(P_{k,l}(s)(x_{\alpha}-x_{\beta}))}{dt^l}\rangle dt=0.\]
So the sequence $\{P_{k,l}(s)(x_{\alpha})\}$ converge and
$P_{k,l}(s)$ is a compact operator.

Now we prove that $P_{m,m}(s)$ is Fredholm and then our lemma is proved. If $p_{m,m}(s,t)$ is positive definite
for each $s,t\in[0,1]$ , we can choose $p_{k,l}(s,t)$ such that $\Ii_{s,R}$ is positive definite for each $s$. So
$P_{m,m}(s)$ is a compact perturbation of a Fredholm operator and is Fredholm. Here it is only required that
$p_{m,m}(s,t)$ is continuous in $t$. In the general case, we have to assume that $p_{m,m}(s,t)$ is $C^m$ in $t$.
Consider the operator $p_{m,m}(s,\cdot):H\to H$. Let $j:H_R\to H$ be the injection. Then $p_{m,m}(s,\cdot)$ is
invertible and $p_{m,m}(s,\cdot)j$ is Fredholm. For any $x\in H_R$ and $y=H$, the inner product $\langle
(P_{m,m}(s)-p_{m,m}(s,\cdot))x,y\rangle_m$ consists only the lower-order terms (i.e., no second-order differential
involved) and some boundary terms. Similar to the above proof, we can conclude that the lower-order terms
correspond to compact operators. The boundary terms correspond to finite rank operators. So
$jP_{m,m}(s)-p_{m,m}(s,\cdot)j$ is compact. Since $p_{m,m}(s)j$ and $j$ are Fredholm, $jP_{m,m}(s)$ and
$P_{m,m}(s)$ are Fredholm. \qe

The following lemma is the key to the proof of Theorem \ref{t05m2.1}.

\begin{lemma} \label{l05m5.2} {\rm (i)} Any solution $u\in H^1([0,T];\C^{2mn})$ of (\ref{e05m2.10})
can be expressed by $u=u_{p_s,x}$ for some $x\in H^m([0,T];\C^n)$,
and the following three conditions are equivalent:
\begin{enumerate}
\item[{\rm (a)}] $x\in\ker\;{\cal I}_{s,R}$;
\item[{\rm (b)}] $x\in\ker L_{s,W_{2m}(R)}$;
\item[{\rm (c)}] $u_{p_s,x}$ is a solution of (\ref{e05m2.10}) and
$(u_{p_s,x}(0),u_{p_s,x}(T))\in W_{2m}(R)$.
\end{enumerate}

{\rm (ii)} If $p_s$ is $C^1$ in $s$, then for any $x,y\in
H^m([0,T];\C^n)$, we have
\begin{equation} \label{e05m5.1}\left\langle\left(\frac{d}{ds}p_s\right)\bar u_{0,x},
\bar u_{0,y}\right\rangle
=-\left\langle\left(\frac{d}{ds}b(p_s)\right)u_{p_s,x},
u_{p_s,y}\right\rangle.
\end{equation}

{\rm (iii)} Let $J\in\GL(\C^m)$ be skew self-adjoint, and
$b_s(t)\in\gl(\C^m)$, $0\le s\le 1$, $0\le t\le T$ is a continuous
family of self-adjoint matrices. Let $\gamma_s$ be the fundamental
solutions of the linear Hamiltonian system
\begin{equation}\label{e05m5.2}-J{\dot u}=b_su.\end{equation}
If $b_s$ is $C^1$ in $s$, we have
\begin{equation}\label{e05m5.3}
\frac{\partial}{\partial t}(-J\gamma_s^{-1}\frac{\partial
\gamma_s}{\partial s})=\gamma_s^*\frac{\partial b_s}{\partial s}
\gamma_s.
\end{equation}

{\rm (iv)} If $p_s$ is $C^1$ in $s$, then for any $x,y\in\ker
L_s$, we have
\begin{equation}\label{e05m5.4}
\left\langle-J_{2m,n}\gamma_{p_s}(T)^{-1}\frac{d\gamma_{p_s}(T)}{d
s}u_{p_s,x}(0),u_{p_s,y}(0)\right\rangle=-\int_0^T
\left\langle\left(\frac{d}{ds}p_s\right)\bar u_{0,x}, \bar
u_{0,y}\right\rangle dt.
\end{equation}
\end{lemma}

\bp (i) The proof for the solution $u$ of (\ref{e05m2.10}) can be
expressed by $u=u_{p_s,x}$ and (a)$\Leftrightarrow$(b) is standard
and we omit it. Now we prove (b)$\Leftrightarrow$(c). By
(\ref{e05m2.6}), we have
$\frac{d}{dt}u_{p_s,x}^k(t)=u_{p_s,x}^{k+1}(t)$ for
$k=0,\ldots,m-2$,
\begin{eqnarray*}\frac{d}{dt}u_{p_s,x}^{m-1}(t)&=&\frac{d^m}{dt^m}x(t)\\
&=&p_{m,m}(s,t)^{-1}u_{p_s,x}^m(t)-\sum_{0\le \beta\le
m-1}p_{m,m}(s,t)^{-1}p_{m,\beta}(s,t)u_{p_s,x}^{\beta}(t)
\end{eqnarray*}
and
\begin{eqnarray*}\frac{d}{dt}u_{p_s,x}^k(t)&=&\sum_{2m-k\le \alpha\le m,0\le\beta\le m}(-1)^{\alpha-m}
\frac{d^{\alpha+k+1-2m}}{dt^{\alpha+k+1-2m}}\left(p_{\alpha,\beta}(s,t)
\frac{d^{\beta}}{dt^{\beta}}x(t)\right)\\
&=&u_{p_s,x}^{k+1}(t)-\sum_{0\le\beta\le m}(-1)^{m+k+1}\left(p_{2m-k-1,\beta}(s,t)\frac{d^{\beta}}{dt^{\beta}}x(t)\right)\\
&=&u_{p_s,x}^{k+1}(t)+(-1)^{m+k}p_{2m-k-1,m}(s,t)p_{m,m}(s,t)^{-1}u_{p_s,x}^m(t)\\
& &+\sum_{0\le\beta\le
m-1}(-1)^{m+k}(p_{2m-k-1,\beta}(s,t)\\
&&-p_{2m-k-1,m}(s,t)p_{m,m}(s,t)^{-1}p_{m,\beta}(s,t))u_{p_s,x}^{\beta}(t)
\end{eqnarray*}
for $k=m,\ldots,2m-1$. Combine the above equations and we get
\begin{equation}\label{e05m5.5}
\frac{d}{dt}u_{p_s,x}(t)=J_{2m,n}b(p_s)u_{p_s,x}(t)+(u_{p_s,x}^{2m}(t),0,\ldots,0).
\end{equation}
By the fact that $L_sx=(-1)^mu_{p_s,x}^{2m}(t)$, we get
(b)$\Leftrightarrow$(c).

(ii) By the definition of $U(p_s)$, $V(p_s)$, $\bar u_{p_s,x}$ and
$\bar u_{0,x}$ in \S\ref{s:main}, direct computation shows
\[V(p_s)^*\left(\frac{d}{ds}p_s\right)V(p_s)=-\frac{d}{ds}P(p_s).\]
Thus for all $x,y\in H_R$, we have
\begin{eqnarray*}
\left\langle\left(\frac{d}{ds}p_s\right)\bar u_{0,x}, \bar
u_{0,y}\right\rangle &=&-\left\langle
U(p_s)^*\left(\frac{d}{ds}P(p_s)\right)U(p_s)\bar u_{0,x},
\bar u_{0,y}\right\rangle \\
&=&-\left\langle \left(\frac{d}{ds}P(p_s)\right)\bar
u_{p_s,x},\bar u_{p_s,y}\right\rangle \\
&=&-\left\langle\left(\frac{d}{ds}b(p_s)\right)u_{p_s,x},
u_{p_s,y}\right\rangle .
\end{eqnarray*}

(iii) By the definition of $\gamma_s$, we have
$\gamma_s^*J\gamma_s=J$, and
\begin{eqnarray*}
\frac{\partial}{\partial
t}(-J\gamma_s^{-1}\frac{\partial\gamma_s}{\partial s})&=&
J\gamma_s^{-1}{\dot
\gamma_s}\gamma_s^{-1}\frac{\partial\gamma_s}{\partial s}
-J\gamma_s^{-1}\frac{\partial^2\gamma_s}{\partial
s\partial t}\\
&=&J\gamma_s^{-1}(-J^{-1}b_s)\frac{\partial\gamma_s}{\partial
s}-J\gamma_s^{-1}\frac{\partial}{\partial s}(-J^{-1}b_s\gamma_s)\\
&=&J\gamma_s^{-1}J^{-1}\frac{\partial b_s}{\partial s}\gamma_s\\
&=&\gamma_s^*\frac{\partial b_s}{\partial s}\gamma_s.
\end{eqnarray*}

(iv) follows from (ii), (iii) and the fact that
$\gamma_{p_s}u_{p_s,x}(0)=u_{p_s,x}$ for all $x\in\ker L_s$. \qe

Now we can prove Theorem \ref{t05m2.1}.

We begin with a simple case.

\begin{lemma}\label{l05m5.3}
Let $\Ii_{Id,R}$ be the inner product on $H_R$. If $\epsilon>0$
satisfies $[-\epsilon,0]\cap\sigma(p_{m,m}(0,t))=\emptyset$ for
all $t\in [0,T]$, we have
\begin{equation}\label{e05m5.6}
-\sf\{\Ii_{0,R}+a\Ii_{Id,R}; a\in [0,\epsilon]\}
=i_{W_{2m}(R)}(\{\gamma_{p_s+aI_{(m+1)n}}(T);0\le a\le T\}).
\end{equation}
\end{lemma}

\bp By Lemma \ref{l05m5.1}, $\Ii_{0,R}+a\Ii_{Id,R}$, $a\in
[0,\epsilon]$ is a continuous family of Fredholm quadratic forms.
By the definition of the spectral flow we have
\begin{equation}\label{e05m5.7}
\sf\{\Ii_{0,R}+a\Ii_{Id,R}; a\in
[0,\epsilon]\}=\sum_{a\in(0,\epsilon]}\dim\ker
(\Ii_{0,R}+a\Ii_{Id,R}).
\end{equation}

Set
\[Z_a=-J_{2m,n}\left(\gamma_{p_s+aI_{(m+1)n}}(T)\right)^{-1}\frac{d\gamma_{p_s+aI_{(m+1)n}}(T)}{d
a}\] for $a\in[0,\epsilon]$. By (iv) of Lemma \ref{l05m5.2}, $Z_a$
is non positive definite. Let $v\in\C^{2mn}$ be a vector such that
$\langle Z_a v,v\rangle=0$. By (i) of Lemma \ref{l05m5.2}, there
exists $x\in\ker L_s$ such that $v=u_{p_s+aI_{(m+1)n},x}(0)$. By
(iv) of Lemma \ref{l05m5.2}, we have $\bar u_{0,x}(t)=0$ for all
$t\in[0,T]$. Thus $x=0$, $u_{p_s+aI_{(m+1)n},x}=0$ and $v=0$. So
$Z_a$ is negative definite. By Lemma \ref{l05m4.3}, Proposition
4.1, (i) of Lemma \ref{l05m5.2} and the definition of Maslov-type
index we have
\begin{eqnarray}
i_{W_{2m}(R)}(\{\gamma_{p_s+aI_{(m+1)n}}(T);0\le a\le T\})
&=&-\sum_{a\in(0,\epsilon]}\dim\Gr\left((\gamma_{p_s+aI_{(m+1)n}}(T))\cap
W_{2m}(R)\right)\nn\\
&=&-\sum_{a\in(0,\epsilon]}\dim\ker
(\Ii_{0,R}+a\Ii_{Id,R}).\label{e05m5.8}
\end{eqnarray}

Combine (\ref{e05m5.7}) and (\ref{e05m5.8}), we get
(\ref{e05m5.6}). \qe

\ss{Proof of Theorem \ref{t05m2.2} and Corollary \ref{c05m2.1}
}\label{ss:pf-t2.2}

We now in the position to prove Theorem \ref{t05m2.1}.

{\bf Proof of Theorem \ref{t05m2.1}.}\hspace{2mm} We divide the
proof into two steps.

{\bf Step 1.} We apply Proposition 4.2. Set
\[A_s=L_s^*,\quad D_m=H_0^{2m}([0,T];\C^n),\quad D_M=H^{2m}([0,T];\C^n).\]
Then $A_s$ is injective for each $s$ and $L_{s,W_{2m}(R)}$, $0\le
s\le 1$ is a continuous family of self-adjoint operators. Define
the trace map $\hat\gamma:D_M\to\C^{4mn}$ by
$\hat\gamma(x)=(u_{p_s,x}(0),u_{p_s,x}(T))$ for $x\in D_M$. Then
$\hat\gamma$ induce an isomorphism $\D_M/D_m\to\C^{4mn}$. After
identify the two space $\D_M/D_m$ and $\C^{4mn}$, we have
$\hat\gamma=\gamma$. Direct computation shows
\[\omega_s(x+D_m,y+D_m)=\langle
J_{2m,n}u_{p_s,x}(0),u_{p_s,y}(0)\rangle-\langle
J_{2m,n}u_{p_s,x}(T),u_{p_s,y}(T)\rangle.\] Let $D_s$ be the
domain of $L_{s,W_{2m}R}$. Then $\gamma(D_s)=W_{2m}(R)$ and
$\gamma(\ker A_s^*)=\Gr(\gamma_{p_s}(T))$. By Proposition 4.2 we
have
\begin{eqnarray}
-\sf\{L_{s,W_{2m}(R)};0\le s\le
1\}&=&\Mas\{W_{2m}(R),\Gr(\gamma_{p_s}(T));0\le s\le
1;\omega_s\}\nn\\
&=&\Mas\{\Gr(\gamma_{p_s}(T)),W_{2m}(R);0\le s\le
1;-\omega_s\}\nn\\
&=&i_{W_{2m}(R)}(\{\gamma_{p_s}(T);0\le s\le 1\}).\label{e05m5.9}
\end{eqnarray}

{\bf Step2.} We claim that
\begin{equation}-\sf\{I_{s,R};0\le s\le
1\}=i_{W_{2m}(R)}(\{\gamma_{p_s}(T);0\le s\le
1\}).\label{e05m5.10}
\end{equation}

Let $\Ii_{Id,R}$ be the inner product on $H_R$. Let $\epsilon>0$
be small enough such that
$[-\epsilon,0]\cap\sigma(p_{m,m}(s,t))=\emptyset$ for all
$(s,t)\in [0,1]\times[0,T]$. By Lemma \ref{l05m5.1}, $\sf\{{\cal
I}_s+a\Ii_{Id,R}\}$ is well-defined. For each $c\in[0,1]$, there
exist $\delta_c>0$ and $\epsilon_c\in(0,\epsilon]$ such that
$\ker({\cal I}_s+\epsilon_c\Ii_{Id,R})=\{0\}$ for all
$s\in(c-\delta_c,c+\delta_c)\cap[0,1]$.

Let $[s_0,s_1]$ be a subinterval of
$(c-\delta_c,c+\delta_c)\cap[0,1]$. Consider the spectral flow
$\sf\{{\cal I}_s+a\Ii_{Id,R}\}$ and the Maslov-type index
$i_{W_{2m}(R)}(\gamma_{p_s+aI_{(m+1)n}}(T))$. Because of the
homotopic invariance of spectral flow and Maslov-type index, both
integers must vanish for the boundary loop going counter clockwise
around the rectangular domain from the corner point $(s_0,0)$ via
the corner points $(s_1,0)$, $(s_1,\epsilon_c)$, and
$(s_0,\epsilon_c)$ back to $(s_0,0)$. The spectral flow and Maslov
index vanish on the top segment of our box. By the preceding
lemma, the left and the right side segments of our curves yield
vanishing sum of spectral flow and Maslov index. So, by the
additivity under catenation, we have \[-\sf\{I_{s,R};s_0\le s\le
s_1\}=i_{W_{2m}(R)}(\{\gamma_{p_s}(T);s_0\le s\le s_1\}).\]

Since $[0,1]$ is compact, there exist $c_0,\ldots,c_{N-1}\in[0,1]$
and a partition of $0=s_0<s_1<\ldots<s_N=1$ of $[0,1]$ such that
$[s_j,s_{j+1}]\subset(c_j-\delta_{c_j},c_j+\delta_{c_j}]$ for
$j=0,\ldots N-1$. Then (\ref{e05m5.10}) follows from additivity
under catenation of spectral flow and Maslov-type index.

{\bf Step 3.} Since $\gamma_{p_s}(0)=I_{2mn}$, by the homotopic
invariance of Maslov-type index we have
\begin{equation}\label{e05m5.11}i_{W_{2m}(R)}(\{\gamma_{p_s}(T);0\le s\le
1\})=i_{W_{2m}(R)}(\gamma_{p_1})-i_{W_{2m}(R)}(\gamma_{p_0}).
\end{equation} \qe

{\bf Proof of Theorem \ref{t05m2.2}.}\hspace{2mm} We divide the
proof into three steps.

{\bf Step 1.} (\ref{e05m2.14}), (\ref{e05m2.15}) holds for $C^1$
path $\gamma$ with $\gamma_0=I_{2m}$.

Set $H=L^2([0,T];\C^n)$ and $H_R=\{x\in H;(x(0),x(T))\in R\}$. Let
$F_R$ be a closed operator on $H$ with domain $H_{R^K}$ defined by
$F_Rx=-K\dot x$ for all $x\in H_R$. Set
\[X=L^2([0,T],\C^{2n}),
\quad D_{W_K(R)}=\{x\in H^1([0,T];\C^{2n});(x(0),x(t))\in
W(R)\}.\] Let $A_{W_K(R)}\in{\cal C}(X)$ with domain $D_{W_K(R)}$
be defined by $A_{W_K(R)}x=-J_K\dot x$ for $x\in D_{W_K(R)}$. Let
$b(t)\in\gl(\C^{2n})$ and $C\in\B(X)$ be defined by
$b(t)=-J_K\dot\gamma(t)\gamma(t)^{-1}$, $t\in[0,T]$ and
$(Cx)(t)=b(t)x(t)$ for $x\in X$, $t\in[0,T]$. Then we have
$F_R^*=-F_{R^K}$.

Consider the standard orthogonal decomposition
\[\C^{2n}=(\C^n\times\{0\})\oplus(\{0\}\times \C^n).\]
It induces orthogonal decompositions $X=H\oplus H$ and
$D_{W_K(R)}=H_{R^K}\oplus H_R$. Under such orthogonal
decompositions, $A_{W_K(R)}$ is in block form
$A_{W_K(R)}=\pmatrix{0&F_R^*\cr F_R&0\cr}$. Let $C$ be in block
form
\[C=\pmatrix{C_{1,1}&C_{1,2}\cr C_{2,1}&C_{2,2}\cr}.\]

By the definition of $b(t)$ and the symplectic path $\gamma$ we
have
\[b(t)=\pmatrix{K^*({\dot M}_{2,1}M_{1,1}^{-1}-{\dot M}_{2,2}M_{2,2}^{-1}M_{2,1}M_{1,1}^{-1})&K^*{\dot M}_{2,2}M_{2,2}^{-1}\cr
-K{\dot M}_{1,1}M_{1,1}^{-1}&0\cr}.\] Since $M_{2,2}^*KM_{1,1}=K$,
we have $K^*{\dot M}_{2,2}M_{2,2}^{-1}=-(M_{1,1}^*)^{-1}{\dot
M}_{1,1}^*K^*$. So there holds
\begin{eqnarray*}
K^*({\dot M}_{2,1}M_{1,1}^{-1}-{\dot
M}_{2,2}M_{2,2}^{-1}M_{2,1}M_{1,1}^{-1})&=&K^*{\dot
M}_{2,1}M_{1,1}^{-1}+(M_{1,1}^*)^{-1}{\dot
M}_{1,1}^*K^*M_{2,1}M_{1,1}^{-1}\\
&=&(M_{1,1}^*)^{-1}\left(\frac{d}{dt}(M_{1,1}^*K^*M_{2,1})\right)M_{1,1}^{-1}.
\end{eqnarray*}

Clearly we have
\[\ker(F_R-C_{2,1})=\{M_{1,1}x(0);(x(0),M_{1,1}(T)x(0))\in
R^K\}.\] Since $\ind(F_R-C_{2,1})=\ind F_R=\dim(\Gr(I_{mn})\cap
R^K)-\dim(\Gr(I_{mn})\cap R)$, we have
\[\dim\ker(F_R-C_{2,1})^*=\dim S(T)+\dim(\Gr(I_{mn})\cap
R)-\dim(\Gr(I_{mn})\cap R^K).\]
Let $x,y\in\ker(F_R-C_{2,1})$.
Then we have
\begin{eqnarray*}
\langle C_{1,1}x,y\rangle
&=&\int_0^T\left\langle(M_{1,1}^*)^{-1}\left(\frac{d}{dt}(M_{1,1}^*K^*M_{2,1})\right)M_{1,1}^{-1}x,y\right\rangle
dt\\
&=&\int_0^T\left\langle(M_{1,1}^*)^{-1}\left(\frac{d}{dt}(M_{1,1}^*K^*M_{2,1})\right)
M_{1,1}^{-1}M_{1,1}x(0),M_{1,1}y(0)\right\rangle
dt\\
&=&\int_0^T\left\langle\left(\frac{d}{dt}(M_{1,1}^*K^*M_{2,1})\right)x(0),y(0)\right\rangle
dt\\
&=&\langle M_{1,1}(T)^*K^*M_{2,1}(T)x(0),y(0)\rangle.
\end{eqnarray*}
By Proposition \ref{p05m3.4}, Proposition \ref{p05m4.3} and the
definition of $S(t)$, we have (\ref{e05m2.14}) and
\begin{eqnarray*}
i_{W_K(R)}(\gamma)&=&-\sf\{A_{W_K(R)}-sC;0\le s\le 1\}\\
&=&m^+((M_{1,1}(T)^*K^*M_{2,1}(T))|_{S(T)})+\dim(\Gr(I_{mn})\cap
R^K)-\dim S(T).
\end{eqnarray*}

{\bf Step 2.} Define the set
\[Y=\{M\in\GL(\C^{2n});M=\pmatrix{M_{1,1}&0\cr
M_{2,1}&M_{2,2}\cr},M^*J_KM=J_K\}.\]
Note that any symplectic loop
$\gamma$ in $Y$ is homotopic to the loop in $Y$ starting from
$I_{2n}$. By the homotopic invariance of the Maslov-type index and
Step 1, we have $i_{W_K(R)}(\gamma)=0$ for any loop in $\gamma$ in
$Y$. For a general $\gamma$ in $Y$, we can connect $I_{2n}$ and
the endpoints $\gamma(0)$ and $\gamma(T)$ in $Y$ by $C^1$ paths.
Then (\ref{e05m2.14}) follows from Step 1 and the path additivity
of Maslov-type index under catenation. \qe

Now we turn to the proof of Corollary \ref{c05m2.1}.

{\bf Proof of Corollary \ref{c05m2.1}.}\hspace{2mm} Let
$x=(x_0,\ldots,x_{m-1})$ and $y=(y_0,\ldots,y_{m-1})$ be two
vectors in $\C^{mn}$. By direct calculation we get our form of
$\gamma_{p_0}=(\gamma_{k,l}(t))_{k,l=0,\ldots,2m-1}$ and
(\ref{e05m2.16}) with $p_{m,m}(0,t)=p_{m,m}(1,t)$. Then we have
\begin{eqnarray*}\langle M_{1,1}(T)^*K_{m,n}^*M_{2,1}(T)x,y\rangle
&=&\sum_{k,l=0,\ldots,m-1}\left\langle\left(\frac{1}{(m-k-1)!(m-l-1)!}\right.\right.\\
& &\left.\left.\int_0^T t^{2m-k-l-2}(p_{m,m}(1,t))^{-1}dt\right)x_l,y_k\right\rangle\\
&=&\int_0^T\left\langle(p_{m,m}(1,t))^{-1}\sum_{l=0,\ldots,m-1}\frac{t^{m-l-1}}{(m-l-1)!}x_l,\right.\\
& &\left.\sum_{k=0,\ldots,m-1}\frac{t^{m-k-1}}{(m-k-1)!}y_k\right\rangle dt\\
\end{eqnarray*}
Since $p_{m,m}(1,t)$ is positive definite for each $t\in[0,T]$, we
have $\langle M_{1,1}(T)^*K_{m,n}^*M_{2,1}(T)x,x\rangle\ge 0$. If
$\langle M_{1,1}(T)^*K_{m,n}^*M_{2,1}(T)x,y\rangle=0$, we have
$\sum_{k=0,\ldots,m-1}\frac{t^{m-k-1}}{(m-k-1)!}x_k=0$ for all
$t\in[0,T]$. By taking derivative with $t$, we have
$\sum_{l=0,\ldots,k}\frac{t^{k-l-1}}{(k-l-1)!}x_l=0$ for all
$k=0,\ldots, m-1$ and $t\in[0,T]$. Then we get $x_k=0$ for
$k=0,\ldots,m-1$ and $x=0$. Thus $M_{1,1}(T)^*K_{m,n}^*M_{2,1}(T)$
is positive definite.

Let $p_s=(1-s)p_0+sp_1$. Clearly $\Ii_{0,R}$ and $L_{0,W_{2m}(R)}$
is non negative definite. For sufficiently large $r>0$, we have
\[\langle L_{s,W_{2m}(R)}x,x\rangle=\Ii_{s,R}(x,x)+r\langle
x,x\rangle>0\] for each $x\ne 0$ in the domain of
$L_{s,W_{2m}(R)}$ and $s\in[0,T]$. Then $L_{s,W_{2m}(R)}+rI$ is
positive definite for each $s\in[0,1]$. Note that
$M_{1,1}(0)=I_{mn}$ and $S(0)=S$. By the definition of the
spectral flow, Theorem \ref{t05m2.1} and Theorem \ref{t05m2.2}, we
have
\begin{eqnarray*}
m^-(\Ii_{1,R})&=&-\sf\{\Ii_{s,R};0\le s\le 1\}\\
&=&i_{W_{2m}(R)}(\gamma_{p_1})-i_{W_{2m}(R)}(\gamma_{p_0})\\
&=&i_{W_{2m}(R)}(\gamma_{p_1})-(\dim S(T)+\dim S(0)-\dim S(T))\\
&=&i_{W_{2m}(R)}(\gamma_{p_1})-\dim S\\
&=&-\sf\{L_{s,W_{2m}(R)};0\le s\le 1\}\\
&=&m^-(L_{s,W_{2m}(R)}). \end{eqnarray*} \qe

\ss{Proof of Theorem \ref{t05m2.3}}\label{ss:pf-t2.3}

Let $a$, $p_1$, $p_1^{'}$ and $R^{'}$ be as in {\S}\ref{s:main}.
Firstly we prove (\ref{e05m2.18}). The following lemma follows
from direct calculation.

\begin{lemma} \label{l05m5.4} We have
\begin{eqnarray}
\label{e05m5.12} p_1^{'}&=&\pmatrix{a^{*}&0\cr {\dot
a}^{*}&a^{*}\cr} p_1
\pmatrix{a&\dot a\cr 0&a\cr},\\
\label{e05m5.13}
b(p_1^{'})&=&\diag(a^{-1},a^{*})b(p_1)\diag(a^{*-1},a)
+\pmatrix{0&-a^{-1}{\dot a}\cr -{\dot a}^{*}a^{*-1}&0\cr}.
\end{eqnarray}
\end{lemma}\qe

By Corollary \ref{c05m4.2} we have

\begin{corollary}\label{c05m5.1} We have
\begin{equation} \label{e05m5.14}
\gamma_1^{'}=\diag(a^*,a^{-1})\gamma_1\diag(a(0)^{*-1},a(0)).
\end{equation}
\end{corollary}






{\bf Proof of Theorem \ref{t05m2.3}.}\hspace{2mm} By the
definition of $R^{'}$ we have
\[(R^{'})^{2,b}=\{(x,y)\in\C^{2n};(a(0)^*x,a(T)^*y)\in R^{2,b}\}.\]
By Theorem \ref{t05m2.2} and
Lemma \ref{l05m4.4}, we have
\begin{eqnarray*} i_{W_2(R^{'})}(\gamma_1^{'}) &=&
i_{W_2(R^{'})}(\diag(a^*,a^{-1})\gamma_1\diag(a(0)^{*-1},a(0)))\\
&=&i_{W_2(R)}(\gamma_1)+i_{W_2(R^{'})}(\diag(a^*,a^{-1})\diag(a(0)^{*-1},a(0)))\\
&=&i_{W_2(R)}(\gamma_1)+\dim (\Gr(I_n)\cap (R^{'})^{2,b}))
-\dim(\Gr(a(T)^*a(0)^{*-1})\cap(R^{'})^{2,b})\\
&=&i_{W_2(R)}(\gamma_1)+\dim (\Gr(I_n)\cap (R^{'})^{2,b}))
-\dim(\Gr(I_n)\cap R^{2,b}).
\end{eqnarray*} \qe





\smallskip

{\bf Acknowledgements.} This work was partially done when the
author visited MIT in 2001 and MPI, Leipzig in 2002-2003. The
research atmosphere is very nice there. The author sincerely thank
Professor Gang Tian and Professor Chun-gen Liu for helpful
discussion and valuable suggestions, and the referees for their
careful reading, pointing out mistakes and typos, and valuable
comments on the earlier version of this paper.
\bibliographystyle{abbrv}

\end{document}